\documentclass[preprint,12pt]{elsarticle}

\usepackage{graphicx,amsmath,amsfonts,color,amssymb}
\usepackage{tikz}
\usepackage{cleveref}

\def\bfm#1{\boldsymbol{#1}} %bold math 
\def\RR{\mathbb{R}}
\def\NN{\mathbb{N}}

\def\ZZ{\mathbb{Z}}

\newcommand{\im}{\operatorname{im}\nolimits}

\newtheorem{thm}{Theorem}
\newtheorem{lem}{Lemma}
\newtheorem{dfn}{Definition}

\newproof{pf}{Proof}
\newproof{pot}{Proof of Theorem \ref{thm2}}
\journal{Computer Aided Geometric Design}

\usetikzlibrary{intersections}
\usetikzlibrary{arrows.meta}

\begin{document}

\begin{frontmatter}

\title{Geometric approximation of the sphere by triangular polynomial spline patches}

\author[address1,address2]{Ale\v{s} Vavpeti\v{c}}
%\author[address1]{Ale\v{s} Vavpeti\v{c}}
\ead{ales.vavpetic@fmf.uni-lj.si}

\author[address1,address2]{Emil \v{Z}agar\corref{corauth}}
%\author[address1]{Emil \v{Z}agar\corref{corauth}}
\ead{emil.zagar@fmf.uni-lj.si}
\cortext[corauth]{Corresponding author}

\address[address1]{Faculty of Mathematics and Physics, University of Ljubljana, 
Jadranska 19, Ljubljana, Slovenia}
\address[address2]{Institute of Mathematics, Physics and Mechanics, Jadranska 19, Ljubljana, Slovenia}

\begin{abstract}
  A sphere is a fundamental geometric object widely used in (computer aided) geometric design.
  It possesses rational parameterizations but no parametric polynomial parameterization exists. 
  The present study provides an approach to the optimal approximation of 
  equilateral spherical triangles by parametric polynomial patches
  if the measure of quality is the (simplified) radial error. As a consequence, optimal approximations 
  of the unit sphere by parametric polynomial spline patches underlying on particular regular spherical
  triangulations arising from a tetrahedron, an octahedron and an icosahedron inscribed in the unit sphere
  are provided. Some low total degree spline patches with corresponding geometric smoothness are 
  analyzed in detail and several numerical examples are shown confirming the quality of approximants. 
\end{abstract}


\end{frontmatter}

%\linenumbers

\section{Introduction}\label{sec:introduction}
It is well known that a sphere possesses a parametric rational parameterization \cite{Schicho-Rational-Surfaces-98}
but no parametric polynomial parameterization exists. Since a sphere is one of fundamental geometric objects and rational parameterizations
might be sometimes difficult to handle, it is worth to study its parametric polynomial approximations. Obviously, it is enough to consider 
the unit sphere, since any other one can be translated
and scaled to the unit one without affecting the
quality of an approximant. An obvious approach is to construct 
a parametric polynomial spline approximation based on parametric polynomial approximations of spherical triangles tiling a (part of a) sphere.
The obtained spline should be at least continuous but its smoothness is also often required. 
Since the approximating spline is usually
used for visualisation, it is enough to consider geometric continuous splines 
($G^k$ continuous splines \cite{FarinHoschekKim-02-Handbook}, \cite{Kiciak-2017-Gk-book}).
In comparison with $C^k$ continuous splines  they offer some aditional free parameters which can be used to improve the quality of 
the approximation or for modelling.\\
There are not many references dealing with this topic available in the literature. 
The pioneering paper on a geometric interpolation of a general surface can be found in \cite{Morken-2005-parametric-surface}.
Its specific type of interpolation was studied in \cite{Jaklic-Kozak-Krajnc-Vitrih-Zagar-2006-patch}. 
An optimal approximation of symmetric surfaces by biquadratic B\'{e}zier surfaces is in 
\cite{Eisele-1994-best-biquadratic}. Recently, a special type of geometric
interpolation by parametric polynomials was presented in \cite{Jaklic-Kanduc-2017-Argyris}. In this study we concentrate on optimal geometric 
approximation of the unit sphere based on geometric interpolation.\\
The paper is organized as follows. In \Cref{sec:preliminaries} some preliminaries are presented. 
A general approach to the parametric polynomial approximation of equilateral spherical triangles is
given in \Cref{sec:equilateral_triangles}. A detailed overview of the geometric continuity of parametric patches is 
provided in \Cref{sec:Gk}. 
In Sections \ref{section:quadraticG0}--\ref{section:quarticG2} particular cases of geometric 
approximation are studied, i.e., a quadratic $G^0$, a cubic $G^1$, a cubic $G^2$, a quartic $G^1$ and a quartic $G^2$. 
The paper is concluded by \Cref{sec:closure}.

\section{Preliminaries}\label{sec:preliminaries}

The main goal of this paper is an optimal approximation of a given sphere by geometric continuous 
parametric polynomial splines composed by triangular parametric polynomial patches. 
As already mentioned before, it is enough to consider the unit sphere
${\mathcal S}$.

\begin{figure}[!htb]
    \centering
    \includegraphics[width=0.35 \linewidth]{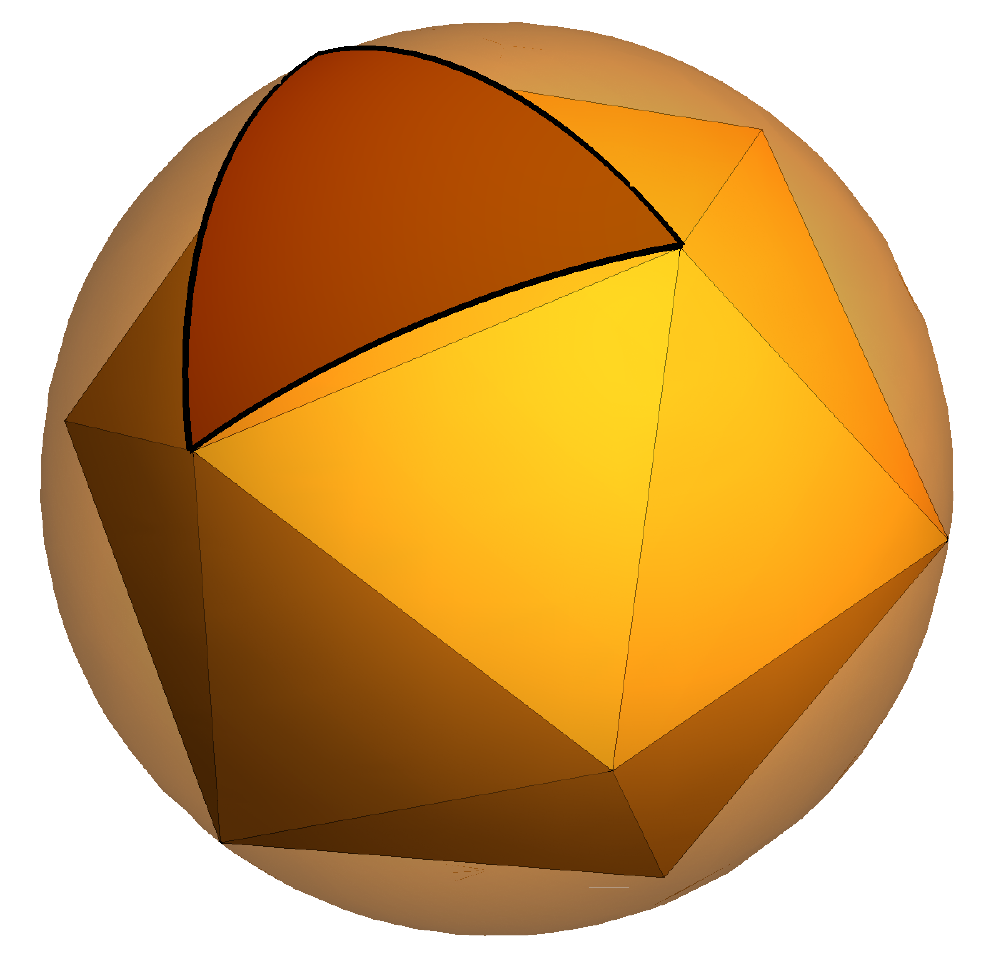}
  \caption{An underlying polyhedron ${\mathcal P}$ (the inner solid), 
  the unit sphere $\mathcal{S}$ (transparent) and a corresponding spherical triangle.}
  \label{fig:spherical_triangulation}
\end{figure}

An approximation of the sphere ${\mathcal S}$ will rely on its particular triangulation by spherical triangles 
related to an underlaying polyhedron. Thus let ${\mathcal P}$ be a polyhedron with the origin 
of the coordinate system in its interior,
vertices $V:=V_{\mathcal P}$ on ${\mathcal S}$ and triangular faces $T:=T_{\mathcal P}$.
The radial projection of triangles from $T$ onto the sphere ${\mathcal S}$ induces  its
triangulation by spherical triangles (see \Cref{fig:spherical_triangulation}).
In order to construct a good (or an optimal) parametric polynomial spline approximant of (a part of) ${\mathcal S}$,
one has to find a good (or an optimal) approximant of each projected triangle from $T$ by triangular parametric polynomial
patch. These patches should be then put together in a smooth parametric polynomial spline patch.
Note that the number of free parameters of a parametric 
polynomial spline patch depends on the cardinality  of $T$, on the degree $n$ of the spline and on the 
order $k$ of smoothness. This makes the problem extremely difficult to be solved in general. 
Thus we shall restrict to some special polyhedra ${\mathcal P}$ inducing spherical triangulations which are built of 
congruent equilateral spherical triangles. 
It is known that in this case ${\mathcal P}$ must be a tetrahedron, an octahedron or an icosahedron
(see, e.g., \cite{Coxeter-1973-Regular-polytopes}). 
Consequently, due to the symmetry reasons, the number of free parameters reduces dramatically and 
this gives some hope that the optimal solution can be found.

\section{Approximation of equilateral spherical triangles}\label{sec:equilateral_triangles}
\label{sec:approximation_of_equilateral_spherical_triangles}
%Let  a spherical cap of a unit sphere be given by the standard parametrization 
%$\bfm{s}:[-\pi,\pi]\times [0,2\pi]:\to \RR^3$, where
%\begin{equation}\label{eq:unit_sphere}
%  \bfm{s}(\psi,\phi)=(\cos\psi\cos\phi,\cos\psi\sin\psi,\sin\phi)^T.
%\end{equation}
  In this section an optimal approximation of an equilateral spherical triangle by a triangular parametric
  polynomial patch will be considered. Without loss of generality we can consider the particular
  equilateral spherical triangle $\bfm{s}_c$ 
  with the mass point at $(0,0,1)$
  %i.e., at the northern pole of the unit sphere centered at the origin, 
  and its vertices given as
   \begin{equation*}\label{eq:equilateral_triangle_coordinates}
     \bfm{v}_0=\left(c,0,\sqrt{1-c^2}\right)^T,
     \quad 
     \bfm{v}_1=\left(-\tfrac{1}{2} c,\tfrac{\sqrt{3}}{2} c,\sqrt{1-c^2}\right)^T,
     \quad  
     \bfm{v}_2=\left(-\tfrac{1}{2} c,-\tfrac{\sqrt{3}}{2} c,\sqrt{1-c^2}\right)^T,
  \end{equation*}
  where $c:=\cos\psi$ and $\psi\in\left[0,\tfrac{\pi}{2}\right)$.
  Note that if $\bfm{s}_c$ is a projection of one of the faces of tehrahedron, octahedron or icosahedron, 
  then the angle $\psi$ implies
 \begin{equation}\label{eq:particular_c}
    c=\frac{2\sqrt{2}}{3},\ c=\frac{\sqrt{6}}{3}\ {\rm and }\ 
    c=\sqrt{\frac{2\left(5-\sqrt{5}\right)}{15}},
  \end{equation}
respectively. 
  Since the image of $\bfm{s}_c$ is a part of the unit sphere,
  it does not possess an exact parametric polynomial parametrization. Thus the goal of this section is to construct 
  an optimal tirangular parametric polynomial approximants of $\bfm{s}_c$.
  It is convenient to consider them in Bernstein-B\'{e}zier representation.
  Let $\Delta$ be the 2-simplex parametrized as 
  $\Delta=\left\{(u,v)\in\RR^2;\ 0\leq u,v\leq 1,\ v\leq 1-u\right\}$.
  A triangular parametric polynomial patch $\bfm{p}:\Delta\to\RR^3$ 
  of total degree $n\in\NN$ can then  
  be defined as
  \begin{equation}\label{eq:parametric_patch}
    \bfm{p}(u,v)=\sum_{i+j+k=n}B_{i,j,k}(u,v)\bfm{b}_{i,j,k},\quad i,j,k\in\ZZ_+,
  \end{equation}
  where
  \begin{equation*}
    B_{i,j,k}(u,v)=\frac{n!}{i! j! k!}u^i v^j(1-u-v)^k
  \end{equation*}
  are bivariate Bernstein polynomials and $\bfm{b}_{i,j,k}\in\RR^3$ are corresponding 
  control points. If one is looking for good triangular parametric polynomial patch approximants 
  \eqref{eq:parametric_patch} of $\bfm{s}_c$, it is 
  natural to require that vertices of a patch coincide with vertices of $\bfm{s}_c$, i.e., 
  \begin{equation}\label{eq:vertices_of_patch}
    \bfm{b}_{n,0,0}=\bfm{v}_0,\quad
    \bfm{b}_{0,n,0}=\bfm{v}_1,\quad
    \bfm{b}_{0,0,n}=\bfm{v}_2. 
  \end{equation}
  If we also require that the control points
  $\bfm{b}_{i,j,0}$ are in the plane passing through the origin, $\bfm{v}_0$ and $\bfm{v}_1$, 
  the control points $\bfm{b}_{0,i,j}$ are 
  in the plane passing through the origin, $\bfm{v}_1$ and $\bfm{v}_2$, and
  the control points $\bfm{b}_{i,0,j}$ are in the plane passing through the origin, $\bfm{v}_0$ and $\bfm{v}_2$, 
  we are able to measure the distance $d_r$ 
  between $\bfm{s}_c$ and $\bfm{p}$ as
  \begin{equation*}\label{def:dr}
    d_r(\bfm{s}_c,\bfm{p}):=\max_{(u,v)\in\Delta}\left|\|\bfm{p}(u,v)\|-1\right|,
  \end{equation*}
  where $\|\cdot\|$ is the Euclidean norm.
  This can be considered as the maximal radial distance between the point on the spherical 
  triangle $\bfm{s}_c$ and the corresponding point on the triangular parametric polynomial patch $\bfm{p}$
  in the radial direction. The function $\|\bfm{p}\|-1$ is an irrational function thus we usually define
  the simplified radial distance
  \begin{equation}\label{def:ds}
    d_s(\bfm{s}_c,\bfm{p}):=\max_{(u,v)\in\Delta}\left|\|\bfm{p}(u,v)\|^2-1\right|,
  \end{equation}
  which is a scalar polynomial of the total degree $2n$. Note that $d_r$ and $d_s$ 
  share the same location of zeros and extrema but, in general, these two measures do not provide the same optimal approximant.\\
  Before we proceed, let us explain the following observation which significantly simplifies some
  technical computations. The proof is easy and will be omitted.
  
\begin{figure}[!htb]
  \minipage{0.49\textwidth}
    \centering
    \includegraphics[width=0.6 \linewidth]{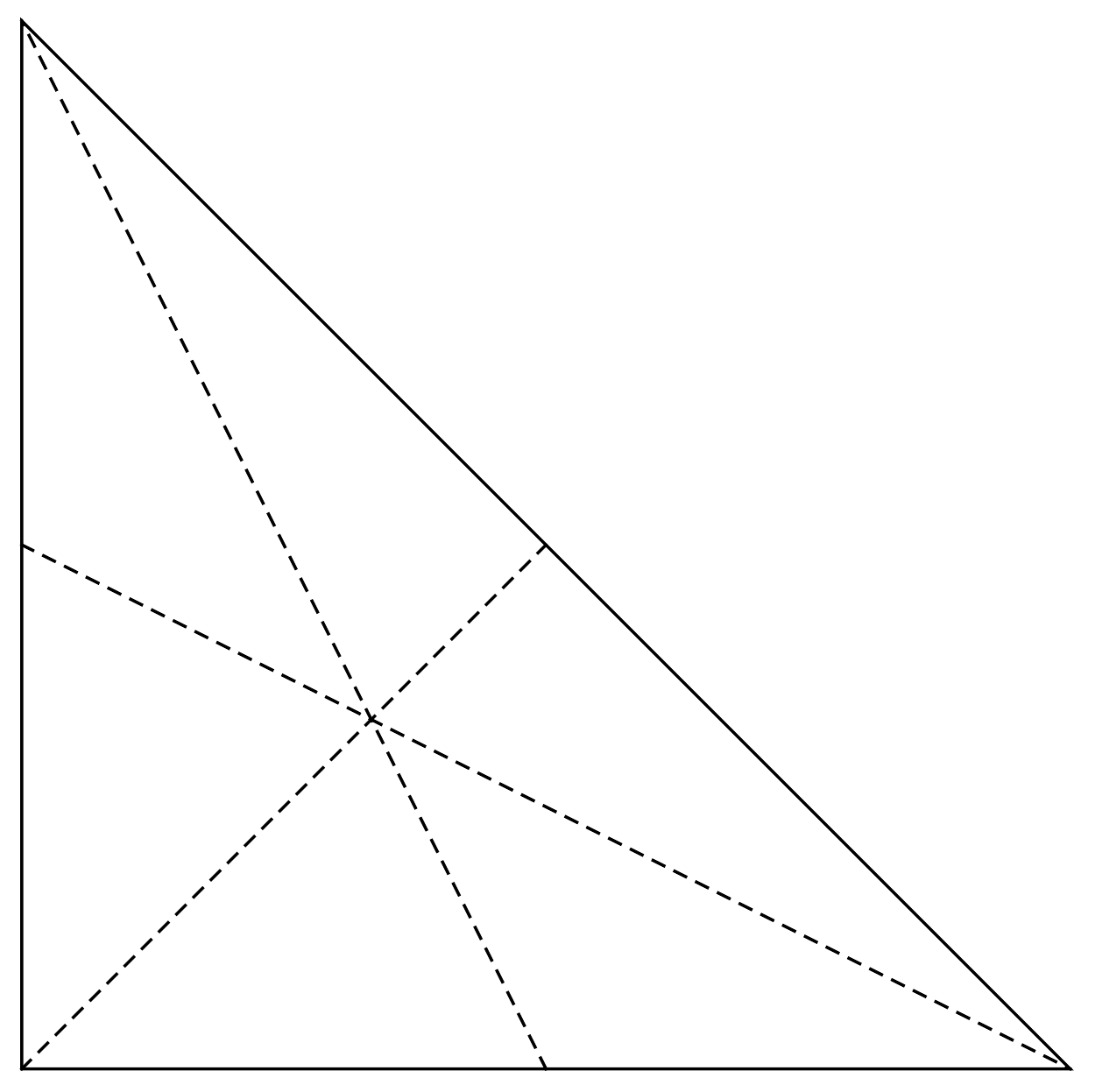}
  \endminipage\hfill
  \minipage{0.49\textwidth}
    \centering
    \includegraphics[width=1\linewidth]{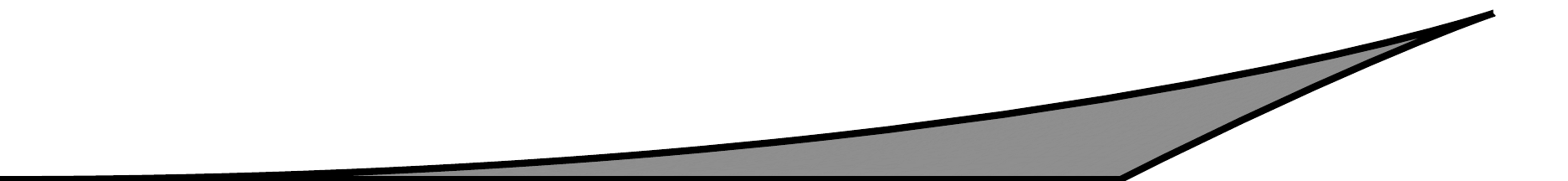}
  \endminipage\hfill
  \caption{The domains $\Delta$ (left) and $\Omega$ (right). Each subtriangle in $\Delta$ is mapped bijectively on
  $\Omega$.}
  \label{fig:delta_omega}
\end{figure}

  \begin{lem}\label{lem:ftog}
    Let $e_1(u,v,w)=u+v+w$, $e_2:=e_2(u,v,w)=uv+uw+vw$ and $e_3:=e_3(u,v,w)=uvw$ be 
    the elementary symmetric polynomials and consider the substitution $w=1-u-v$. If a bivariate
    polynomial $f$ can be written as $f(u,v)=h(e_2,e_3)$, where $h$ is also a bivariate polynomial, then
    $$
     \max_{(u,v)\in\Delta}f(u,v)=\max_{(e_2,e_3)\in \Omega}h(e_2,e_3),
    $$
    where 
    $$\Omega=\Big\{(e_2,e_3)\in[0,\infty)^2;\ 
    \frac{1}{27} \left(9 e_2-2-(2-6 e_2)\sqrt{1-3 e_2}\right) 
    \leq e_3\leq 
    \frac{1}{27} \left(9 e_2-2+(2-6 e_2)\sqrt{1-3 e_2}\right)\Big\}.
    $$
    Moreover, the preimage of the boundary of $\Omega$ under the map 
    $\Delta\to\Omega$, $(u,v)\mapsto (e_2,e_3)$
    is the union of the boundary of $\Delta$ and its medians. Consequently, $h$ has its extrema on the
    boundary of $\Omega$ if and only if $f$ has its extrema on the boundary of $\Delta$ or on its medians.
  \end{lem}
  The domain $\Delta$ consists of six triangles determined by the medians and the sides of $\Delta$,
  and each of them if mapped 
  bijectively on $\Omega$ by the map described in \Cref{lem:ftog} (see \Cref{fig:delta_omega}).
Note that the degree of the polynomial $h$ is significantly smaller than the degree of $f$ and thus the 
analysis of its extrema is much easier which will be of great help later.
  
\section{Geometric continuity of triangular parametric polynomial spline patches}\label{sec:Gk}

Triangular parametric polynomial patches $\bfm{p}$ introduced in the previous section should be composed together
in order to get a triangular parametric polynomial spline patch approximation of (a part of) the unit sphere ${\mathcal S}$.
We are interested in geometrically (or visually) smooth approximation, which is formally defined as follows.
\begin{dfn}\label{def:Gk}
  A parametric patch $\bfm{p}:D\to\RR^3$, where $D\subset \RR^2$ is an open set, is geometric continuous of order $k\in\NN$
  ($G^k$ continuous) if for each point $\bfm{x}\in D$ there exists an open neighbourhood $D_{\bfm{x}}\subseteq D$ of $\bfm{x}$
  and a homeomorphism (a reparameterization)  $\bfm{\rho}_{\bfm{x}}:B^2\to D_{\bfm{x}}$, where $B^2$ is an open ball in $\RR^2$, 
  such that $\bfm{p}\circ\bfm{\rho}_{\bfm{x}}:B^2\to\RR^3$ is regular and $C^k$ continuous.
\end{dfn}
In practice it is not straightforward to check the $G^k$ continuity of a given spline patch. 
Let $\Delta_i$ be the intersection of the unit ball in $\|\cdot\|_1$
norm with the $i$th  quadrant in $\RR^2$ for $i=1,2$. Note that $\Delta=\Delta_1$. 
Let $N_\epsilon(A)=\{\bfm{x}\in\RR^2;\ \min_{\bfm{a}\in A}\|\bfm{x}-\bfm{a}\|_1<\epsilon\}$ 
be the open $\epsilon$ neighbourhood of $A\subset\RR^2$. We are interested in $G^k$ continuity of
a spline patch $\bfm{p}\colon \Delta_1\cup\Delta_2\to \RR^3$ defined by two $C^k$ triangular parametric patches 
$\bfm{p}_i \colon \Delta_i\to\RR^3$, $i=1,2$. 
Note that $C^k$ continuity implies
that there exists an $\epsilon>0$, such that $\bfm{p}_i \colon N_\epsilon(\Delta_i) \to \RR^3$ is $C^k$
continuous for $i=1,2$. Furthermore, a spline patch $\bfm{p}$ is considered as a map defined on
$N_\epsilon(\Delta_1\cup\Delta_2)$ with 
$$\bfm{p}(u,v)=
  \begin{cases}
    \bfm{p}_1(u,v),&u\geq 0,\\
    \bfm{p}_2(u,v),&u < 0.
  \end{cases}
$$
Since we shall consider only continuous spline patches, we can define a $C^k$ continuous 
boundary common curve 
$\bfm{c} \colon (-\epsilon,1+\epsilon) \to \RR^3$ as $\bfm{c}(\tau)=\bfm{p}_1(0,\tau)=\bfm{p}_2(0,\tau)$.

\begin{thm}[\cite{Prautzsch-Boehm-Paluszny-2002}]\label{thm:Boehm-Prautysch-Paluszny}
  Let $\bfm{p}_i \colon N_\epsilon(\Delta_i) \to \RR^3$ be a $C^k$ continuous parametric patch with $k\geq 1$,
  for $i=1,2$. The parametric spline patch $\bfm{p} \colon \Delta_1\cup\Delta_2 \to \RR^3$ is $G^k$ continuous 
  if and only if for every $\tau\in(-\epsilon,1+\epsilon)$ there exists 
  $\delta>0$ and a $C^k$ continuous curve $\bfm{\gamma}_{\tau}:(-\delta,\delta)\to \bfm{p}(N_\epsilon(\Delta_1\cup\Delta_2))$,
  such that
  \begin{itemize}
    \item $\bfm{\gamma}_{\tau}(0)=\bfm{p}(0,\tau)=\bfm{c}(\tau)$,
    \item $\bfm{\gamma}_\tau'(0) \nparallel \bfm{c}'(\tau)$.
  \end{itemize}
\end{thm}
\begin{pf}
  The idea of the proof is the following (see \Cref{fig:Boehm-Prautysch-Paluszny} and \cite[p. 190]{Prautzsch-Boehm-Paluszny-2002}
  for details).
\begin{figure}[!htb]
\begin{center}
\begin{tikzpicture}
\def\a{3}
\def\ep{0.5}
\draw[dotted] (\a+\ep,0)--(0,\a+\ep)--(-\a-\ep,0)--(-\a,-\ep)--(\a,-\ep)--cycle;
\fill[black!20!white] (\a,0)--(0,\a)--(-\a,0)--cycle;
\draw[thick] (0,\a)--(0,0);
\draw[thick] (\a,0)--(0,\a)--(-\a,0)--cycle;
\node at (-\a+\ep,0.4*\ep) {$\Delta_2$};
\node at (\a-\ep,0.4*\ep) {$\Delta_1$};
\node at (-0.9*\a,-1.5*\ep) {$N_\epsilon(\Delta_1\cup \Delta_2)$};
\draw[black] (-0.9,1.8)--(0,1.8)--(1,1.3);
\draw[black] (-0.95,0.7)--(0,0)--(1.4,-0.3);
\node at (-1,0.9) {${\bfm p}^{-1}\circ\bfm{\gamma}_0$};
\node at (0.33*\a,0.35*\a) {${\bfm p}^{-1}\circ\bfm{\gamma}_\tau$};
\begin{scope}[xshift=4cm]
\node at(0,4){$(-\delta,\delta)$};
\draw[-{Latex}] (0.1,3.5)--node[above]{$\bfm{\gamma}_\tau$}++(2,-1.5);
\draw[-{Latex}] (-0.1,3.5)--++(-2,-1.5);
\draw[-{Latex}] (-2,1.8)--node[above]{$\bfm p$}++(4,0);
\end{scope}
\begin{scope}[xshift=9cm]
\draw[dotted] (\a+\ep,0.5*\ep) .. controls +(0,1.8) and +(2.9,0)..(-0.9*\a,\a+0.8*\ep)
.. controls +(0.9,-1.2) and +(0.5,1)..(-\a-1.2*\ep,-0.5*\ep)
 .. controls +(0,-0.2) and +(-0.2,0)..(-\a+\ep,-0.8*\ep)
 .. controls +(1,0.8) and +(-1,0.8)..(\a-0.5*\ep,0)--cycle;
\path[name path=al2,thick] (-0.8*\a+\ep,\a)--(1,0);
\path[name path=al1,thick] (-\a,0) .. controls +(1,0.8) and +(-1,0.8)..(\a,\ep);
\fill[black!20!white] (\a,\ep) .. controls +(0,1) and +(2,0)..(-0.8*\a+\ep,\a)
.. controls +(0.3,-1.2) and +(0.5,1)..(-\a,0)
 .. controls +(1,0.8) and +(-1,0.8)..(\a,\ep);
\draw[thick] (\a,\ep) .. controls +(0,1) and +(2,0)..(-0.8*\a+\ep,\a)
.. controls +(0.3,-1.2) and +(0.5,1)..(-\a,0)
 .. controls +(1,0.8) and +(-1,0.8)..(\a,\ep);
\path[name intersections={of=al1 and al2,by={a}}];
\draw[thick] (a) .. controls +(-0.3,0.8) and +(1,-0.3).. (-0.8*\a+\ep,\a);
\draw [black] plot [smooth] coordinates {(-0.5*\a,0.8*\a) (-0.4*\a,0.8*\a) (-0.1*\a,0.7*\a) (0*\a,0.75*\a)};
\draw [black] plot [smooth] coordinates {(-1.3,1.5) (-0.8,1.1) (a) (1.5,0.6)};
\node at (0.05*\a,0.75*\a) {$\bfm{\gamma}_\tau$};
\node at (-1.55,1.5) {$\bfm{\gamma}_0$};
\end{scope}
\end{tikzpicture}
\end{center}
\caption{A path $\bfm{\gamma}_\tau$ is a composition of an inclusion of the interval $(-\delta,\delta)$
into $N_\varepsilon(\Delta_1\cup\Delta_2)$ and $\bfm p$.
Even though the map 
$\bfm{\gamma}_\tau$ is $G^k$ continuous the corresponding inclusion is usually not even $G^1$ continuous.}
\label{fig:Boehm-Prautysch-Paluszny}
\end{figure}
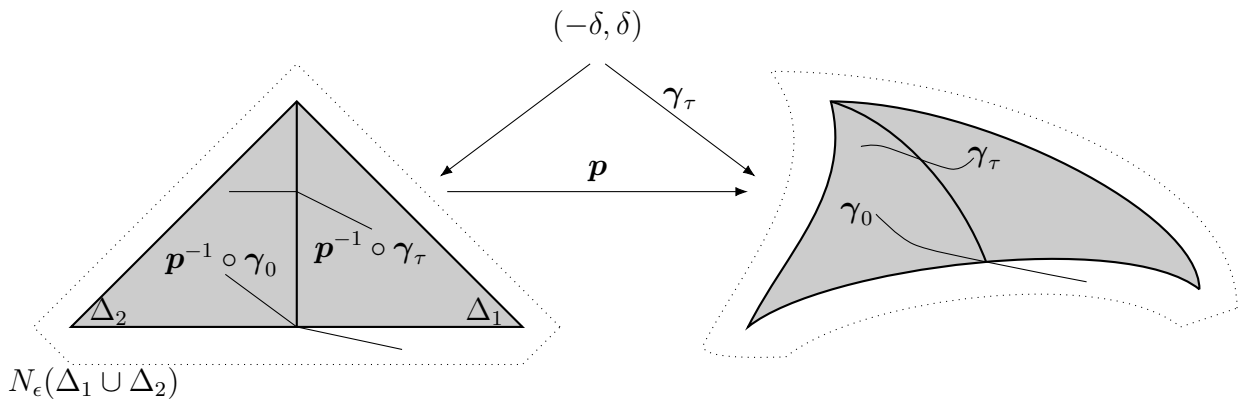
  Let $\bfm{P}=\bfm{\gamma}_\tau(0)=\bfm{c}(\tau)$, $\tau\in (-\epsilon,1+\epsilon)$ 
  be a point on the common curve $\bfm{c}$ of the patches $\bfm{p}_i$. 
  Since they are $G^1$ continuous each of them possesses the tangent plane at $\bfm{P}$.
  But the existence of $\bfm{\gamma}_\tau$ implies that these tangent planes must coincide. 
  Let us denote the common tangent plane by 
  $T_{\bfm{P}}$ and let $\pi_{\bfm{P}}\colon \RR^3\to T_{\bfm{P}}$ be the orthogonal projection. 
  By the implicit function theorem there exists a neighbourhood $U\subset T_{\bfm{P}}$ of $\bfm{P}$ 
  and functions $\bfm{f}_i\colon U\to \bfm{p}_i(N_\epsilon(\Delta_i))$, such that 
  $\pi_{\bfm{P}}\circ\bfm{f}_i=id$, i.e., the image of a patch $\bfm{p}_i(N_\epsilon(\Delta_i))$ 
  is locally a graph of the function $\bfm{f}_i$.
  There exists a neighbourhood $V_1\subset \Delta_1\cap [0,\infty)\times\RR$ and a neighbourhood 
  $V_2\subset \Delta_2\cap (-\infty,0]\times\RR$ of $(0,\tau)$ such that for 
  $U_i=\pi_{\bfm{P}}\circ \bfm{p}_i(V_i)$ we have $U_1\cup U_2=U$ and 
  $U_1\cap U_2=U\cap \pi_{\bfm{P}}(\bfm{c}(-\epsilon,1+\epsilon))$.
  An existence of suitable paths $\bfm{\gamma}_\tau$ implies that the $C^k$ functions ${\bfm{f}_i}|_{U_i}$ induce the $C^k$ function on $U$.
  \qed
\end{pf}

\begin{figure}[!htb]
    \centering
    \includegraphics[width=0.35 \linewidth]{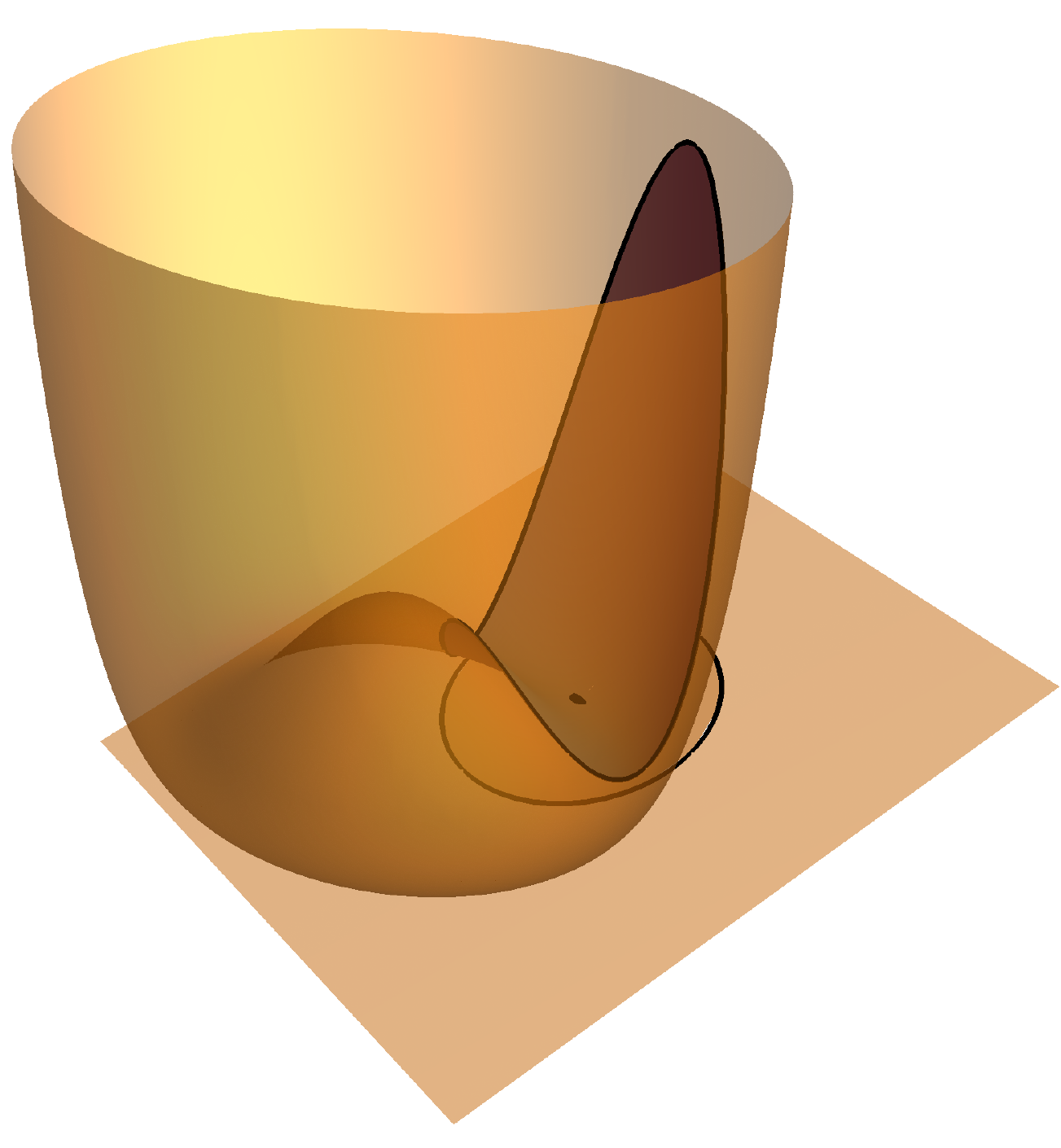}
  \caption{A $G^1$ continuous parametric surface is locally a graph of a function over the tangent plane.}
  \label{fig:graph_over_tangent_plane}
\end{figure}

By the above characterization the $G^k$ continuity of a surface is checked
by confirming the $G^k$ continuity of particular spatial curves which is well understood topic. 
In particular, it is easy to prove that $G^1$ continuity of 
a parametric spline patch $\bfm{p}$ is equivalent to the existence of the tangent plane $\bfm{T}_{\bfm{v}}$ 
at every point $\bfm{v}\in\bfm{p}(N_\epsilon(\Delta_1\cup\Delta_2))$. 
This basically means that for every point
$\bfm{v}\in \bfm{p}(N_\epsilon(\Delta_1\cup\Delta_2))$ the patch $\bfm{p}$ can 
be locally seen as a graph of $C^1$ continuous function $\bfm{f}_{\bfm{v}}$ in the neighbourhood of $\bfm{v}$
in $\bfm{T}_{\bfm{v}}$ (see \Cref{fig:graph_over_tangent_plane}).
Note that once the $G^1$ continuity is established the existence of the tangent plane at every point of the patch follows.
Therefore a spline patch is additionally the $G^k$ continuous, $k\geq 2$, if and only if all functions $\bfm{f}_{\bfm{v}}$
are $C^k$ continuous.

\begin{figure}[!htb]
  \minipage{0.49\textwidth}
    \centering
    \includegraphics[width=0.75 \linewidth]{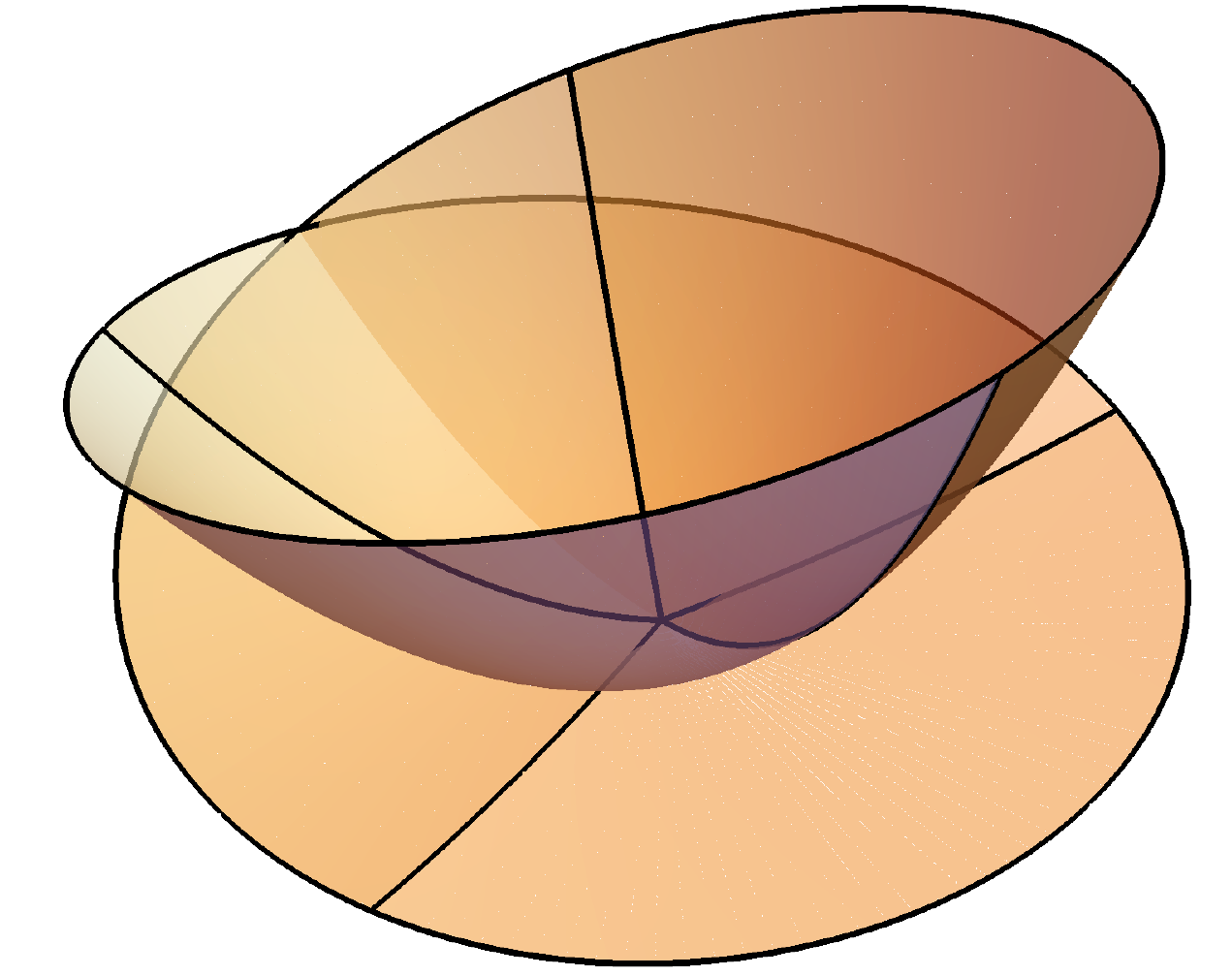}
  \endminipage\hfill
  \minipage{0.49\textwidth}
    \centering
    \includegraphics[width=0.75\linewidth]{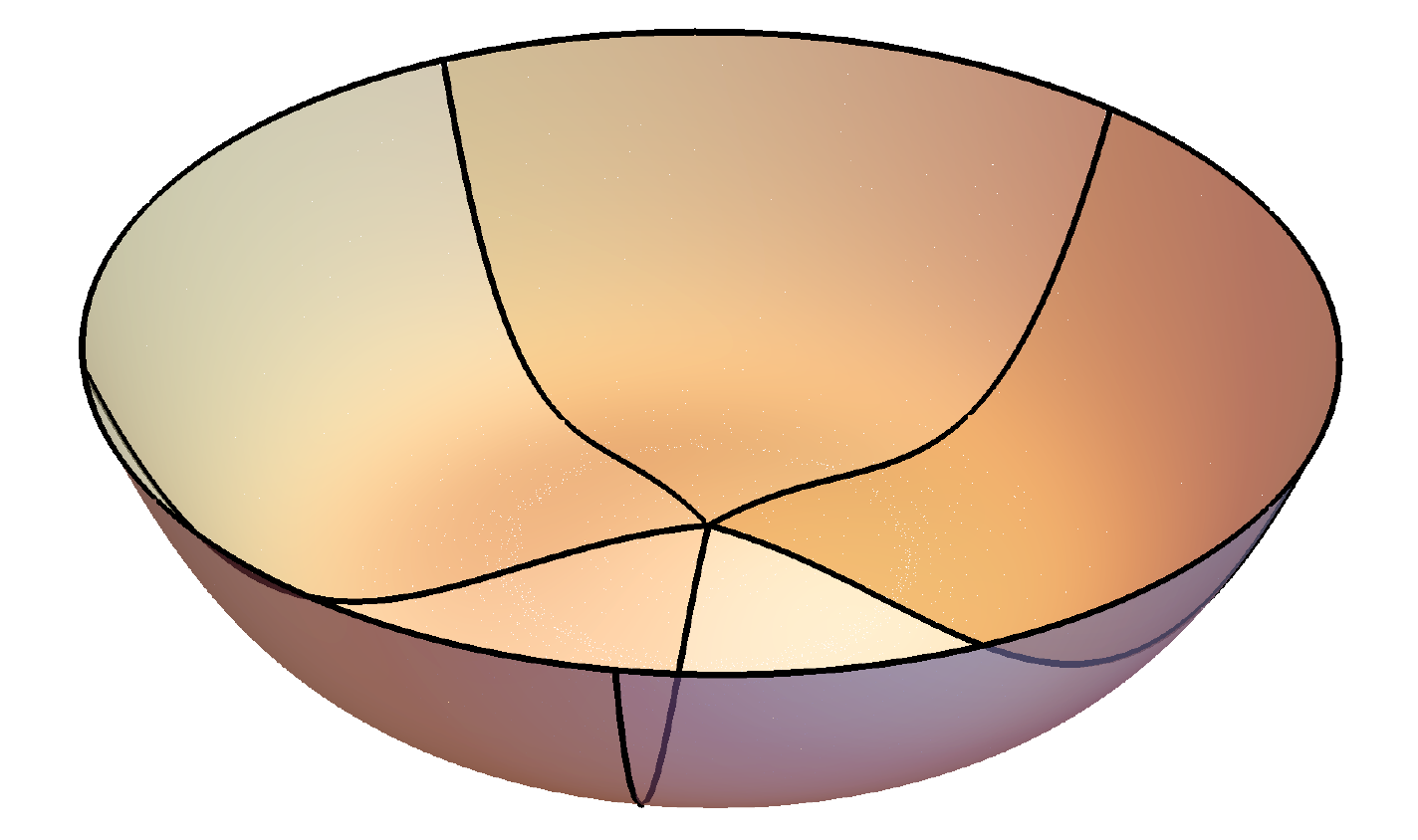}
  \endminipage\hfill
  \caption{A degenerate case of $G^1$ continuous spline patch around an interior point (left),
  and a regular case (right).}
  \label{fig:patch_arround_point}
\end{figure}

We have characterized the $G^k$ continuity of a spline patch composed by two patches sharing a common boundary curve.
When several patches are put together around a common interior point of a continuous spline patch, the $G^k$ continuity must be treated
more carefully (see \Cref{fig:patch_arround_point}).\\ 
Let $\bfm{p}_i \colon N_\epsilon(\Delta) \to \RR^3$, $i=1,2,\dots,r$, be a set of $C^k$ continuous parametric patches
such that $\bfm{p}_i(0,t)=\bfm{p}_{i+1}(t,0)$ for all $t\in(-\varepsilon,1+\varepsilon)$ and all $i=1,\ldots,r$, where $\bfm{p}_{r+1}=\bfm{p}_1$.
Let $\bfm{p}$ be a spline patch induced by patches $\bfm{p}_i$, $i=1,2,\dots,r$.
Let $\bfm{c}_i\colon (-\epsilon,1+\epsilon)\to\RR^2$ be a path defined by $\bfm{c}_i(t)=\bfm{p}_i(0,t)=\bfm{p}_{i+1}(t,0)$.
For every $i$ we can define a continuous patch $\bfm{q}_i\colon \Delta_1\cup\Delta_2\to\RR^3$ such that $\bfm{q}_i|_{\Delta_1}=\bfm{p}_i$ and $\bfm{q}_i|_{\Delta_2}=\bfm{p}_{i+1}\circ R$, where $R\colon\RR^2\to\RR^2$ is the rotation around $(0,0)$ for $-\tfrac\pi 2$.
If all $\bfm{q}_i$ are $G^k$ continuous, $k\ge 1$, all their tangent planes at $\bfm{P}=\bfm{p}_i(0,0)$ coincide.
 We have seen before that there exists a neighbourhood $U$ of $\bfm{P}$ in the resulting tangent plane
 and functions $\bfm{f}_i\colon U\to \RR^3$ such that $\bfm{f}_i$ and $\bfm{f}_{i+1}$ induce a function on $U$ such that $\bfm{q}_i$ is locally its graph. 
It can happen that we cannot find a neighbourhood $V_i\subset \Delta\cap [0,\infty)^2$ such that for $U_i=\pi_{\bfm{P}}\circ \bfm{p}_i(V_i)$ we have $U_i\cap U_{i+1}={\bfm c}_i([0,1+\epsilon))\cap U$, $U_i\cap U_j=\bfm{P}$ for $|i-j|>2$, and $\cup_{i=1}^r U_i=U$
(see \Cref{fig:patch_arround_point}). But if such sets exist then functions $\bfm{f}_i|_{U_i}$ induce a $C^k$ function such that $\im(\bfm{p})$ is locally its graph.

Let us now use the above facts to construct a geometric approximation of two neighbouring 
equilateral spherical triangles.
If $\bfm{p}_1$ is a triangular parametric polynomial patch approximating $\bfm{s}_c$ 
and $\bfm{p}_2=R\bfm{p}_1$ where $R$ is a reflection over one of the boundary B\'ezier curves $\bfm{c}$
of $\bfm{p}_1$, then the spline patch $\bfm{p}_{12}$ 
composed of $\bfm{p}_1$ and $\bfm{p}_2$ is a spline approximation
of two neighbouring spherical triangles similar to $\bfm{s}_c$. Since $R$ preserves $\bfm{c}$
the spline patch $\bfm{p}_{12}$ is automatically $G^0$ continuous. Moreover, by the corollary of
\Cref{thm:Boehm-Prautysch-Paluszny} it is  also $G^1$ continuous if the tangent plane of $\bfm{p}_{1}$
coincides with the tangent plane of $\bfm{p}_{2}$ at each point on $\bfm{c}$. Since $\bfm{p}_{12}$ 
is symmetric over $\bfm{c}$ it is thus enough to verify
that the cross product of the tangent vector of $\bfm{c}$ at any point $\bfm{P}$ on $\bfm{c}$
and the tangent vector of $\bfm{p}_{12}$ in any other direction at 
$\bfm{P}$ is parallel to the radial vector of $\bfm{P}$. Consequently, the obtained cross 
product must be perpendicular to the normal of the plane involving the common boundary B\'ezier curve $\bfm{c}$ 
and the origin. Checking the $G^2$ continuity of the spline patch $\bfm{p}_{12}$ is more challenging and
it will be explained later by using \Cref{thm:Boehm-Prautysch-Paluszny}.\\
If one is looking for a good approximation of the whole unit sphere, the geometric continuity at the 
vertices $V$ where several triangular parametric polynomial patches meet has to be considered. 
Again, the $G^0$ continuity follows from the fact that all relevant patches share the same vertex. 
The $G^1$ continuity is induced by the $G^1$ continuity of two consecutive neighbouring patches since they already determine the common normal of the tangent plane of the spline patch at the considered vertex. 
The $G^2$ continuity also follows from the $G^2$ continuity of two consecutive neighbouring 
patches, since they already determine the two main Gaussian curvatures together with principal directions
of the spline patch at the considered vertex.

Let us now study the geometric continuity of some low degree triangular parametric polynomial 
spline patches. We shall consider only such spline patches which are induced by a 
triangulation of (a part of) a unit sphere 
consisting of equilateral spherical triangles which arise from one of three Platonic polyhedra with all 
its faces congruent to the same equilateral triangle, i.e., tetrahedron, octahedron or icosahedron
(see \cite{Coxeter-1973-Regular-polytopes}).

\section{Quadratic $G^0$ approximation}\label{section:quadraticG0}
  
  Since a linear approximation of the spherical triangle is uniquely defined by the plane passing through its vertices, the first nontrivial case is a quadratic $G^0$ approximation.
  Here $n=2$ and the vertices of the triangular parametric polynomial patch $\bfm{p}$ are given as in   
  \eqref{eq:vertices_of_patch}.
  Thus we are left with three free control points 
  \begin{equation}\label{eq:middle_control_points_quadratic}
    \bfm{b}_{1,1,0},\quad
    \bfm{b}_{0,1,1},\quad 
    \bfm{b}_{1,0,1}.
  \end{equation}
  The spherical triangle $\bfm{s}_c$ possesses several types of symmetry, so must the patch $\bfm{p}$. In particular, $\bfm{s}_c$
  is symmetric with respect to its circular medians, 
  so each control point \eqref{eq:middle_control_points_quadratic} must be on the line through 
  the origin and the midpoint of the line connecting
  corresponding vertices $\bfm{v}_i$, $\bfm{v}_j$, $i\neq j$. Moreover, due to the
  convexity of $\bfm{s}_c$, the control points \eqref{eq:middle_control_points_quadratic}
  must be of the form
  \begin{equation*}\label{eq:optimal_middle_points_quadratic}
    \bfm{b}_{1,1,0}=\frac{\alpha}{2}(\bfm{v}_0+\bfm{v}_1),\quad 
    \bfm{b}_{0,1,1}=\frac{\alpha}{2}(\bfm{v}_1+\bfm{v}_2),\quad 
    \bfm{b}_{1,0,1}=\frac{\alpha}{2}(\bfm{v}_0+\bfm{v}_2),
  \end{equation*}
  where $\alpha\geq 1$ is a free parameter. Note that for $0\leq \alpha<1$ a parametric patch $\bfm{p}$ is inside of the
  underlying tetrahedron and thus it can not be an optimal approximation of the spherical triangle.
  Since $\bfm{p}$ depends on only one unknown parameter $\alpha$, we shall write
  $\bfm{p}(\cdot,\cdot,\alpha):=\bfm{p}$. In order to construct the optimal triangular parametric 
  $G^0$ quadratic approximant of $\bfm{s}_c$ according to the distance \eqref{def:ds}, one has to 
  solve the minimax problem
  $$
    \min_{\alpha>1}d_s(\bfm{s}_c,\bfm{p}(\cdot,\cdot,\alpha)).
  $$
 Let us first observe the error functions $f$ and $g$ defined by $f(u,v,\alpha):=\|\bfm{p}(u,v,\alpha)\|_2^2-1$ and
 $g(u,v,\alpha):=\|\bfm{p}(u,v,\alpha)\|_2-1$. Note that $f$ and $g$ depend also on $c$, 
 but we have fixed it and we shall omit writing $c$ as an argument. Note also that 
 $g=\sqrt{f+1}-1$ thus $f(\cdot,\cdot,\alpha)$ and $g(\cdot,\cdot,\alpha)$ share the same locations of 
 zeros and extrema, and $f(u,v,\cdot)$ and $g(u,v,\cdot)$ share the same monotonicity properties.
 Thus it is enough to consider the properties of $f$ only.
 
\begin{figure}[!htb]
  \minipage{0.33\textwidth}
    \centering
    \includegraphics[width=0.8\linewidth]{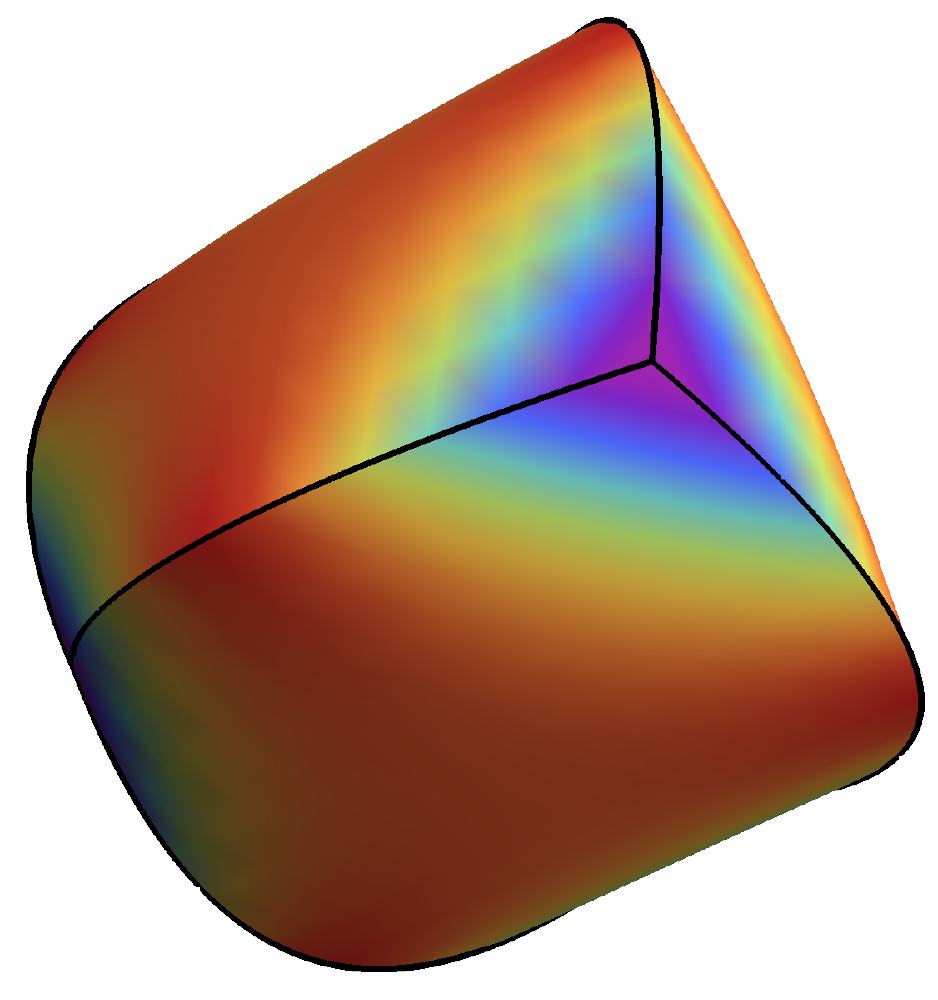}
  \endminipage\hfill
  \minipage{0.33\textwidth}
    \centering
    \includegraphics[width=0.8\linewidth]{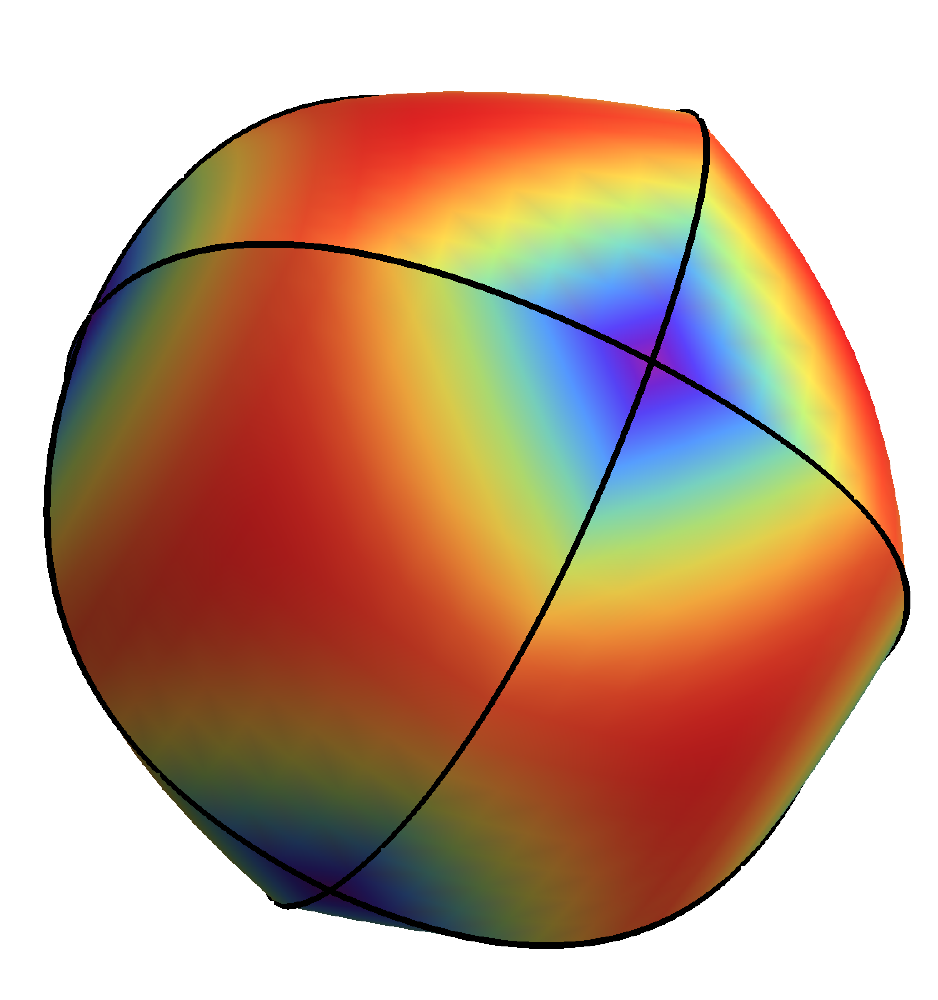}
  \endminipage\hfill
  \minipage{0.33\textwidth}
    \centering
    \includegraphics[width=0.8\linewidth]{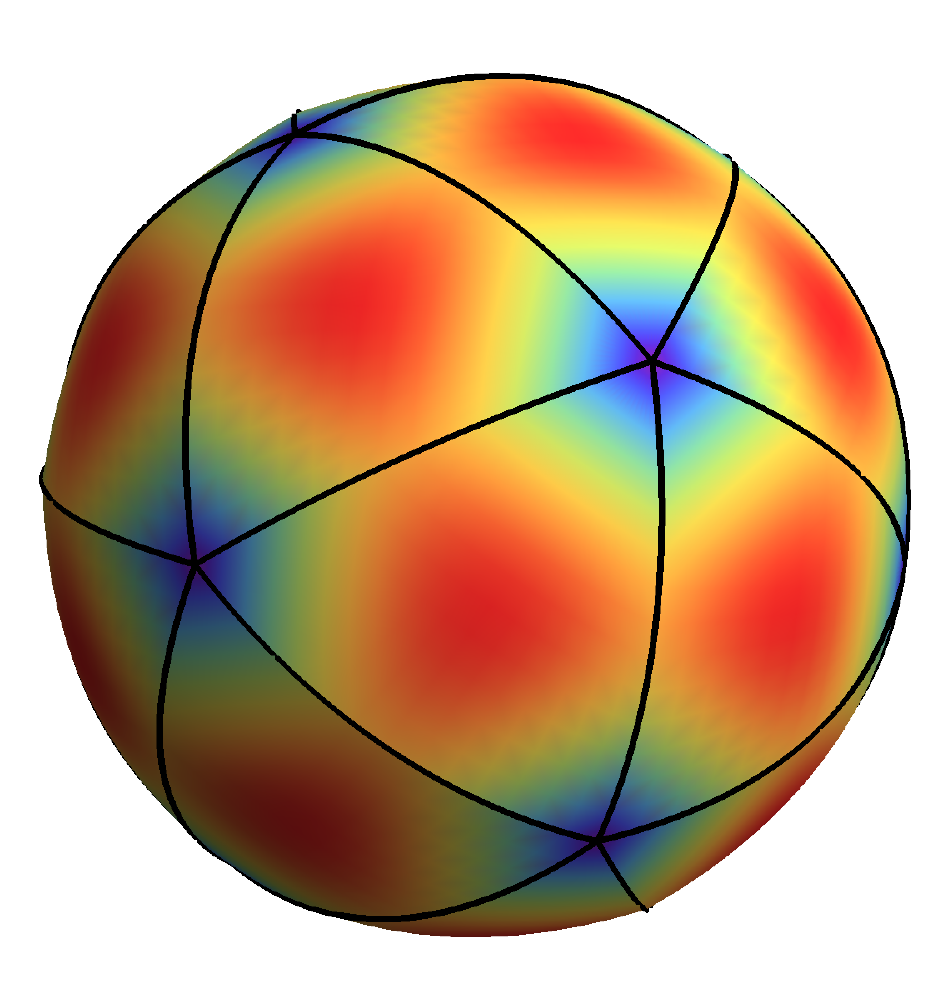}
  \endminipage\hfill
  
  \caption{Figures of the best quadratic $G^0$ approximants of the whole sphere 
  based on the underlying tetrahedron, octahedron and icosahedron (from left to right)
  together with the Gaussian curvatures (red regions indicate higher curvature).
  \label{fig:best_quad_G0}}
\end{figure}

 In order to simplify the analysis, let us define standard symmetric polynomials 
 $e_1(u,v,w)=u+v+w$, $e_2:=e_2(u,v,w)=uv+uw+vw$ and $e_3:=e_3(u,v,w)=uvw$. If we write $w=1-u-v$,
 then it is easy to check that
%  \begin{align*}
%    &f(u,v,w,d)=\left(12 \left(1-c^2\right) e_1(u,v,w)
%      e_3(u,v,w)+\left(4-3 c^2\right) \left(e_2(u,v,w)^2\right.\right.\\
%    &\left.\left.-3 e_1(u,v,w)e_3(u,v,w)\right)\right) d^2\\
%    &+\left(\left(4-3 c^2\right) \left(e_1(u,v,w)^2
%     e_2(u,v,w)-2 e_2(u,v,w)^2-3 e_3(u,v,w) e_1(u,v,w)\right)\right.\\
%    &\left.+12 \left(1-c^2\right) e_1(u,v,w)
%     e_3(u,v,w)\right) d\\
%    &+6 c^2 e_3(u,v,w) e_1(u,v,w)+e_1(u,v,w)^4
%     -4 e_1(u,v,w)^2 e_2(u,v,w)\\
%    &+\left(4-3 c^2\right)e_2(u,v,w)^2-1.
%  \end{align*}
  \begin{align*}
    &f=\left(12 \left(1-c^2\right)
      e_3+\left(4-3 c^2\right) \left(e_2^2-3e_3\right)\right) \alpha^2
    +\left(\left(4-3 c^2\right) 
    \left(e_2-2 e_2^2-3 e_3\right)+12 \left(1-c^2\right) e_3\right) \alpha\\
    &+6 c^2 e_3
     -4 e_2+\left(4-3 c^2\right)e_2^2.
  \end{align*}
  We first observe the following result.
  \begin{lem}\label{lem:f(d)}
    For a fixed point $(u,v)\in\Delta$, and $0<c\leq 1$, the error function 
    $f(u,v,\cdot)$ is an increasing quadratic function on $[0,\infty)$.
  \end{lem}
  \begin{pf}
    To shorten the notation, let us skip writing arguments of $f$. Observe that $f$ can be written as
    $f=\sum_{i=0}^2f_i \alpha^i$
    where
    \begin{align*}
      f_0&=6 c^2 e_3-4e_2+\left(4-3 c^2\right)e_2^2,\\
      f_1&=\left(4-3 c^2\right) \left(e_2-2 e_2^2-3 e_3\right)+12 \left(1-c^2\right) e_3,\\
      f_2&=12 \left(1-c^2\right) e_3+\left(4-3 c^2\right) \left(e_2^2-3e_3\right).
    \end{align*}
    The result of the lemma follows since  $f_2 > 0$ and $ -\tfrac{f_1}{2f_2}\leq 0$ which can be verified
    by a straightforward computations. \qed
  \end{pf}
  Due to the symmetry of the triangular parametric polynomial patch $\bfm{p}$, the necessary and sufficient 
  condition for
  $\max_{(u,v)\in\Delta,\alpha>1}\left|f(u,v,\alpha)\right|$ being minimal is that
  \begin{equation}\label{quad_minimax}
    \min_{(u,v)\in\Delta,\alpha>1}f(u,v,\alpha)=-\max_{(u,v)\in\Delta,\alpha>1}f(u,v,\alpha).
  \end{equation}
  \Cref{lem:f(d)} is crucial for the construction of the best triangular parametric approximant.
  Since $f(u,v,\cdot)$ is an increasing function, there is at most one $\alpha$ for which 
  \eqref{quad_minimax} is fulfilled.
  Thus by \eqref{quad_minimax}, by \Cref{lem:f(d)} and for the reason of symmetry
  we may guess that the parameter $\alpha$ for the optimal
  triangular parametric $G^0$ quadratic approximant is determined by the relation
  $f\left(\tfrac{1}{3},\tfrac{1}{3},\alpha\right)=-f(\tfrac{1}{2},\tfrac{1}{2},\alpha)$ for the simplified radial error and
  by $g\left(\tfrac{1}{3},\tfrac{1}{3},\alpha\right)=-g(\tfrac{1}{2},\tfrac{1}{2},\alpha)$ for the radial error.
  This leads to the following admissible solutions 
  \begin{equation}\label{eq:optimal_d_quadratic}
    \alpha_f=\frac{68-59 c^2-12 \sqrt{196-175 c^2-3 c^4}}{91 c^2-100},\quad 
    \alpha_g=\frac{24-3\sqrt{4-3c^2}-4\sqrt{1-c^2}}{3\sqrt{4-3c^2}+8\sqrt{1-c^2}}.
  \end{equation}
  One should formally check that the obtained solutions really imply the global maxima  
  $f(\tfrac{1}{3},\tfrac{1}{3},\alpha_f)$ and $g(\tfrac{1}{3},\tfrac{1}{3},\alpha_g)$ and the global minima
  $f(\tfrac{1}{2},\tfrac{1}{2},\alpha_f)$ and $g(\tfrac{1}{2},\tfrac{1}{2},\alpha_g)$ over $\Delta$. 
  But this can be easily done by using \Cref{lem:ftog}.
  Optimal parameters together with the corresponding radial distances and minimal and maximal Gaussian
  curvatures for the underlying 
  tetrahedron, octahedron and icosahedron are collected in \Cref{tab:G02}.\\
  The obtained optimal patches $\bfm{p}$ can be used to construct $G^0$ continuous spline patch approximations
  of the sphere defined over triangulation defined by the underlying 
  tetrahedron, octahedron or icosahedron. We just compose $\bfm{p}$ by appropriate rotations.\\
  Let us conclude this section by showing that there is no parameter $\alpha$ for which the above construction would
  give a triangular parametric $G^1$ quadratic spline patch.
  The necessary condition for $G^1$ continuity is that the normal of the tangent plane 
  at the boundary control point, say $\bfm{b}_{2,0,0}$, is parallel to the radius vector of this point. 
  Some easy calculations reveal that this can happen if and only if 
  $$
   \alpha\left((4-3c^2)\alpha-4\right)=0,
  $$
  i.e., if $\alpha=0$ or $\alpha=\tfrac{4}{4-3c^2}$. The first solution is not admissible
  and the second one could, by \eqref{eq:optimal_d_quadratic}, provide the optimal triangular 
  parametric quadratic $G^1$ approximant only for $c=\tfrac{6}{\sqrt{43}}$, which is not one of the considered parameters
  from \eqref{eq:particular_c}. Consequently, there is no optimal triangular parametric
  $G^1$ quadratic spline approximation of the (part of the) unit sphere induced by its
  equilateral spherical triangulation. The optimal $G^0$ spline patches together with their Gaussian
  curvatures are shown on \Cref{fig:best_quad_G0}.

\renewcommand{\arraystretch}{1.5}
\begin{table}[htb]
  \begin{equation*}
    \begin{array}{|l|r|r|r|r|r|}\hline
        & \multicolumn{1}{|c|}{\alpha_f} 
        & \multicolumn{1}{|c|}{\alpha_g} 
        & \multicolumn{1}{|c|}{d_r}
        &\multicolumn{1}{|c|}{K_{min}} 
        & \multicolumn{1}{|c|}{K_{max}}\\ \hline
        \text{tetrahedron} 
        & 3.060496 
        & 3.132163
        & 0.192853
        &-0.19&0.16\\ \hline
        \text{octahedron}
        & 1.965622 
        & 1.968975
        & 0.049691
        &0.01&0.41\\ \hline
      \text{icosahedron} & 1.371294 & 1.371371& 0.008604&0.30&0.71\\ \hline
    \end{array}
  \end{equation*}
  \caption{Optimal parameters $\alpha_f$ according to the simplified radial error, $\alpha_g$ according 
  to the radial error, radial distances $d_r$ according to $\alpha_g$ and the corresponding
  minimal and maximal Gaussian curvatures $K_{min}$, $K_{max}$ for the optimal triangular parametric
  $G^0$ quadratic patches.} \label{tab:G02}
\end{table}

\section{Cubic $G^1$ approximation}\label{section:cubicG1}

In order to improve the quality of the approximation we have to rise the degree of the
triangular parametric polynomial patch. The next practically important case are cubic patches.
Since by \eqref{eq:parametric_patch} the number of control points  rises to $10$ this implies
some new parameters. Assuming \eqref{eq:vertices_of_patch} again, we are left with $7$ 
control points. Due to the symmetries we observe that
\begin{align*}
    \bfm{b}_{2,1,0}&=\alpha\ \bfm{b}_{3,0,0}+\beta\ \bfm{b}_{0,3,0},\quad 
    \bfm{b}_{2,0,1} =\alpha\ \bfm{b}_{3,0,0}+\beta\ \bfm{b}_{0,0,3},\\
    \bfm{b}_{1,2,0}&=\alpha\ \bfm{b}_{0,3,0}+\beta\ \bfm{b}_{3,0,0},\quad 
    \bfm{b}_{0,2,1}=\alpha\ \bfm{b}_{0,3,0}+\beta\ \bfm{b}_{0,0,3},\\
    \bfm{b}_{1,0,2}&=\alpha\ \bfm{b}_{0,0,3}+\beta\ \bfm{b}_{3,0,0},\quad 
    \bfm{b}_{0,1,2}=\alpha\ \bfm{b}_{0,0,3}+\beta\ \bfm{b}_{0,3,0},\\
    \bfm{b}_{1,1,1}&=\left(0,0,\gamma\right)^T.
\end{align*}
The $G^0$ approximation would require the analysis of three-parametric family of triangular
parametric cubic pathches. This is quite a challenging issue, thus we shall focus on
$G^1$ approximation which obviously implies some relations between parameters $\alpha$, $\beta$ and
$\gamma$. Namely, due to the $G^1$ condition at the particular vertex \eqref{eq:vertices_of_patch}, 
say $\bfm{b}_{0,0,3}$, we observe that all directional derivatives at this point must be perpendicular
to its radius vector. This implies $2 \alpha +2 \beta -3 c^2 \beta -2=0$
and 
\begin{equation*}
  \alpha=\frac{1}{2} \left(2-2 \beta +3 c^2 \beta \right).
\end{equation*}
Furthermore, for each $u\in[0,1]$
the normal of the tangent plane of $\bfm{p}$ at $\bfm{p}(u,0)$ must be in the plane
determined by $\bfm{b}_{0,0,3}$, $\bfm{b}_{3,0,0}$ and the origin. This implies the following triples of parameters
\begin{align}
    \alpha_1&=\frac{8-3c^2}{3(4-3c^2)}, &\beta_1&=\frac{4}{3(4-3c^2)},\quad &\gamma_1&=\frac{\sqrt{1-c^2}(8-8c^2+3c^4)}{(1-c^2)(4-3c^2)},\nonumber\\
    \alpha_2&=1, &\beta_2&=0,\quad &\gamma_2&=\frac{\sqrt{1-c^2}(4-3c^2)}{4(1-c^2)},\label{eq:best_cubic_G1}\\
    \alpha_3&=\frac{3c^2}{4-3c^2}, &\beta_3&=\frac{4}{4-3c^2},\nonumber\quad &\gamma_3&=\frac{\sqrt{1-c^2}(8-9c^4)}{(1-c^2)(4-3c^2)}.
\end{align}

It is easy to see that choices $(\alpha_k,\beta_k,\gamma_k)$, $k=2,3$,
imply nonregular patch. For the first one the normal vector vanishes at $u=v=0$, and for the second one the same happens at
$u=v=\tfrac{1}{2}$. The choice $(\alpha_1,\beta_1,\gamma_1)$ provides the only and therefore the best
regular $G^1$ cubic spline patch.
In \Cref{tab:G13} optimal parameters $\alpha_1$, $\beta_1$ and $\gamma_1$, together with the corresponding 
radial distances and minimal and maximal Gaussian curvatures for the underlying 
tetrahedron, octahedron and icosahedron are collected.
%Using \Cref{lem:ftog} it is easy to verify that its maximal 
%simplified radial error 
%is attained at $(\tfrac{1}{3},\tfrac{1}{3})$, namely
%\begin{equation*}
%   f\left(\tfrac{1}{3},\tfrac{1}{3}\right)=\frac{24 c^4-39 c^6+16 c^8}{9 \left(4-3 c^2\right)^2 \left(1-c^2\right)}.
%\end{equation*}

\renewcommand{\arraystretch}{1.5}
\begin{table}[h]
  \begin{equation*}
    \begin{array}{|l|r|r|r|r|r|r|}\hline
        & \multicolumn{1}{|c|}{\alpha_1} 
        & \multicolumn{1}{|c|}{\beta_1}
        & \multicolumn{1}{|c|}{\gamma_1}
        & \multicolumn{1}{|c|}{d_r}
        & \multicolumn{1}{|c|}{K_{min}} 
        & \multicolumn{1}{|c|}{K_{max}}\\ \hline
      \text{tetrahedron} 
      & 1.333333
      & 1.000000
      & 3.666667
      & 0.370370
      & 0.11
      & 3.24\\ \hline
      \text{octahedron} 
      & 1.000000
      & 0.666667
      & 1.732051
      & 0.090551
      & 0.25
      & 1.69\\ \hline
      \text{icosahedron} 
      & 0.793989
      & 0.460655
      & 1.186755
      & 0.016690
      & 0.52
      & 1.25\\ \hline
    \end{array}
  \end{equation*}
  \caption{Optimal parameters $\alpha_1$, $\beta_1$ and $\gamma_1$,
  radial distances $d_r$ and the corresponding
  minimal and maximal Gaussian curvatures $K_{min}$, $K_{max}$ for the optimal triangular parametric
  $G^1$ cubic patches.} \label{tab:G13}
\end{table}

\begin{figure}[!htb]
  \minipage{0.33\textwidth}
    \centering
    \includegraphics[width=0.8 \linewidth]{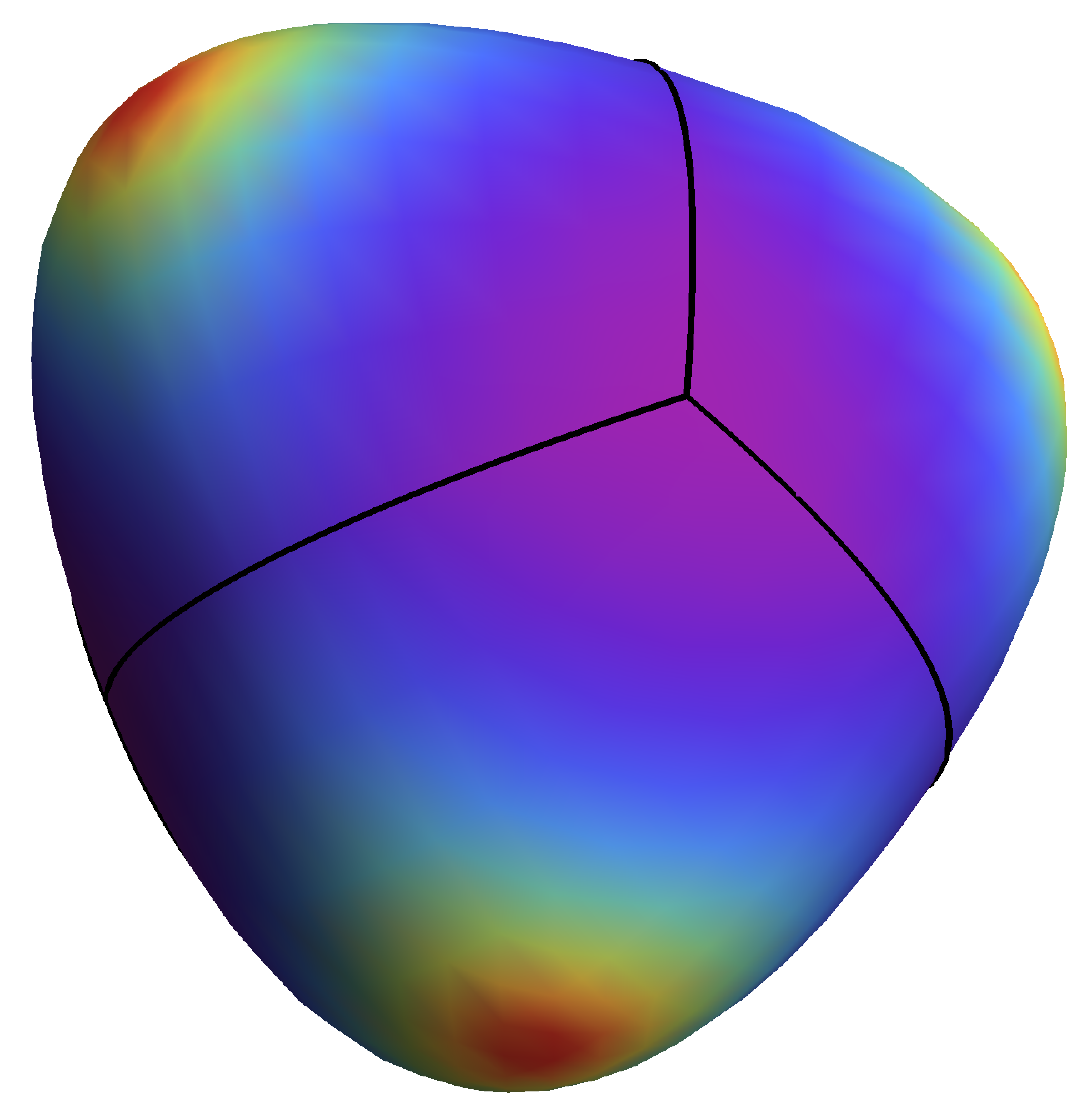}
  \endminipage\hfill
  \minipage{0.33\textwidth}
    \centering
    \includegraphics[width=0.8\linewidth]{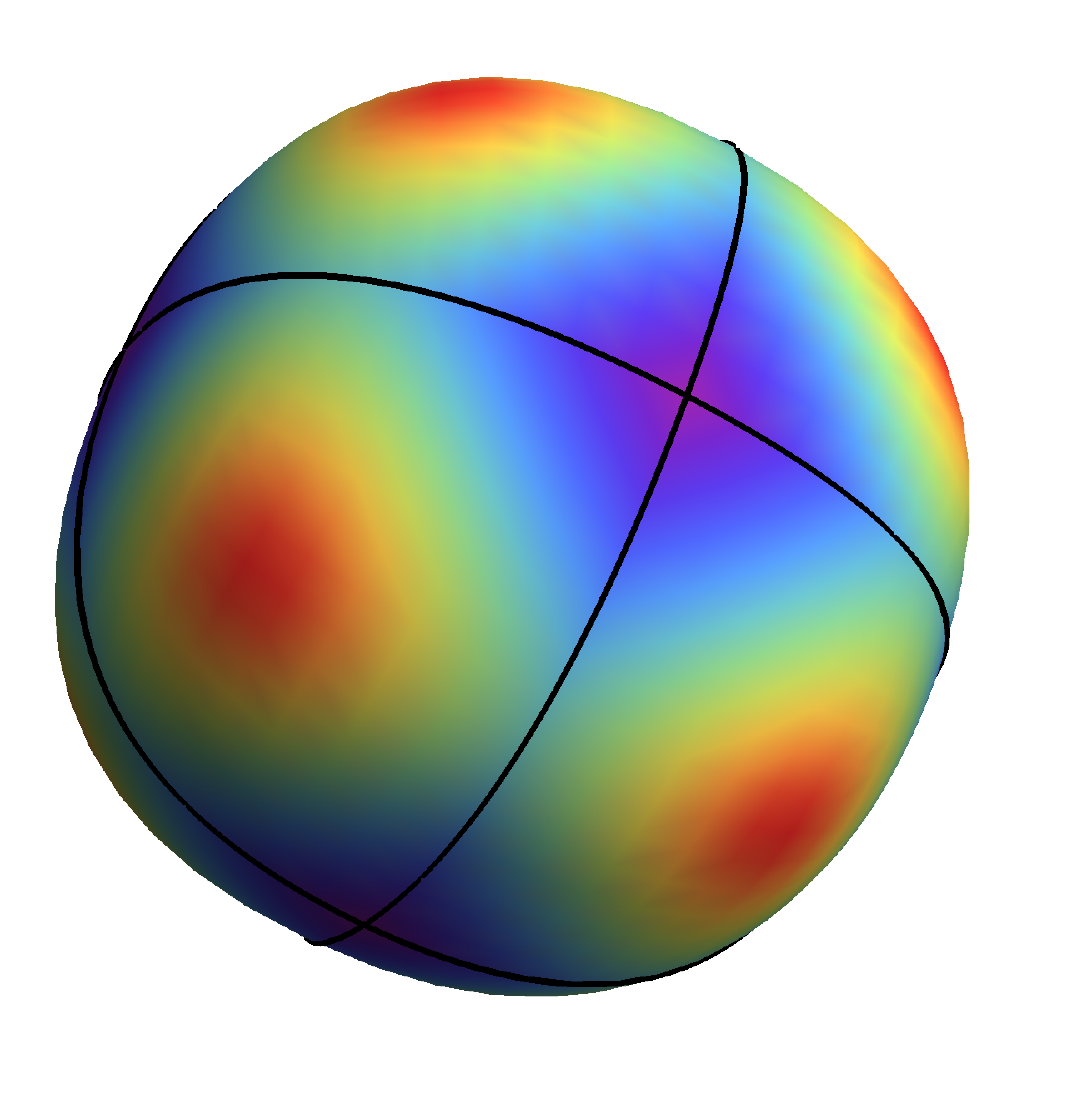}
  \endminipage\hfill
  \minipage{0.33\textwidth}
    \centering
    \includegraphics[width=0.8\linewidth]{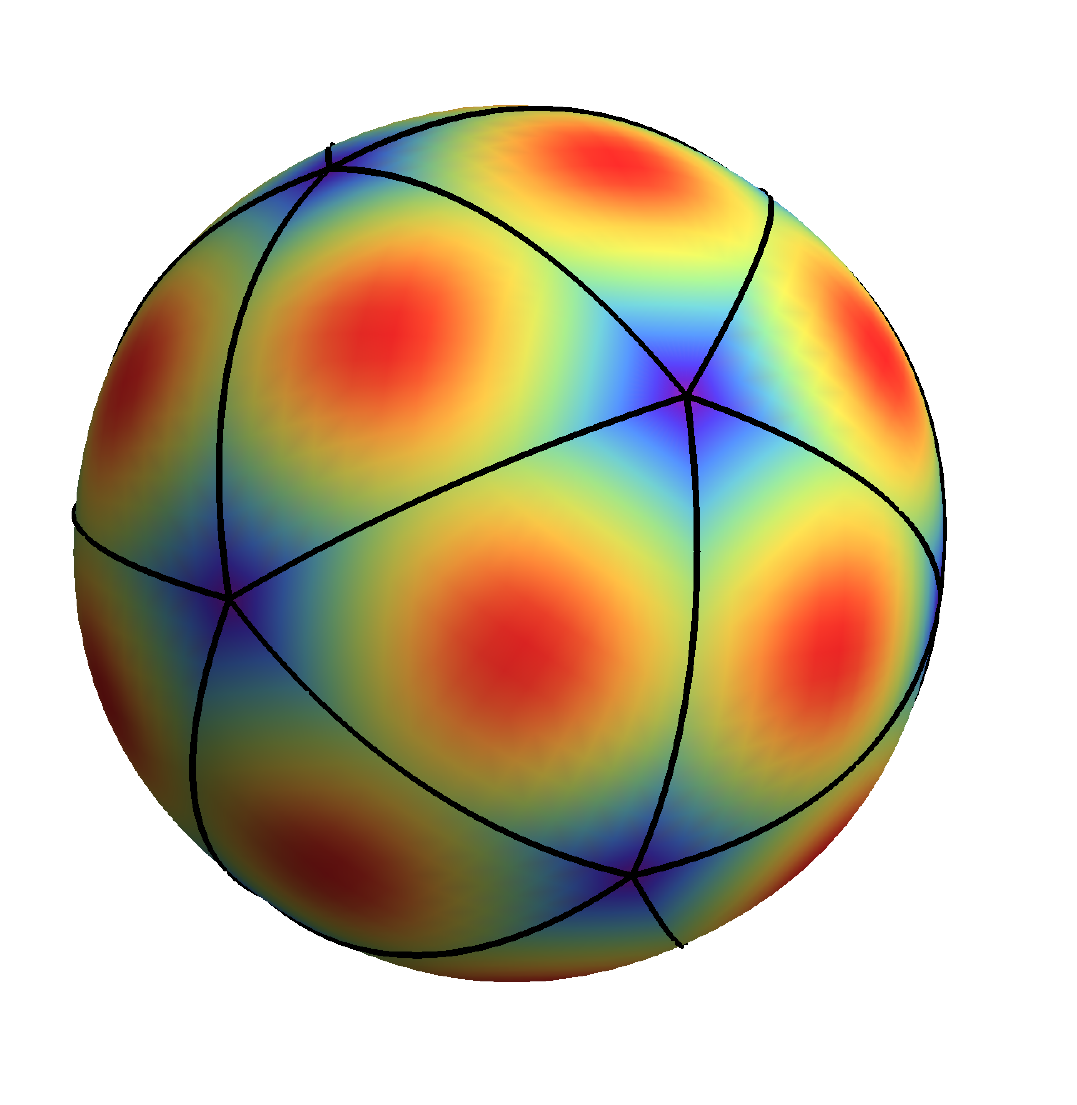}
  \endminipage\hfill
  
  \caption{Figures of the best cubic $G^1$ approximants of the whole sphere 
  based on the underlying tetrahedron, octahedron and icosahedron (from left to right)
  together with the Gaussian curvatures (red regions indicate higer curvature).
  \label{fig:best_cub_G1}}
\end{figure}

The optimal $G^1$ cubic spline patches together with their Gaussian curvatures are shown on \Cref{fig:best_cub_G1}.

\section{Cubic $G^2$ approximation}\label{section:cubicG2}

We have seen in the previous section that for a fixed $c$ the optimal $G^1$ continuous triangular parametric cubic 
spline patches exists with the parameters given by \eqref{eq:best_cubic_G1}.
In the following we shall prove that these spline patches are actually all $G^2$ continuous and consequently
they must be optimal $G^2$ continuous triangular parametric cubic spline approximants of (a part of) the
unit sphere. Since confirming the $G^2$ continuity is much more complicated as checking the
$G^1$ continuity, the underlying theory of $G^2$ continuity will be given in the following subsection.
Moreover, the conditions for $G^k$ continuity for any $k\geq1$ will be derived.

\subsection{$G^k$ continuity of parametric spline patches}\label{subsec:Gk}

In this section the idea how to use \Cref{thm:Boehm-Prautysch-Paluszny} to prove that two
triangular parametric polynomial patches with a common boundary curve form a $G^k$ continuous 
triangular parametric polynomial spline patch will be explained.\\
Let us suppose that two regular triangular parametric polynomial patches $\bfm{p}_i:\Delta\to\RR^3$, $i=1,2$, 
share a common regular $C^k$ continuous boundary curve $\bfm{c}:[0,1]\to\RR^3$. We can assume that
$\bfm{c}=\bfm{p}_1(0,\cdot)=\bfm{p}_2(0,\cdot)$. Let us fix a point $\bfm{c}(\tau)$, $\tau\in[0,1]$.
Our goal is to check the existence of a parametric 
polynomial curve $\bfm{\gamma}:[-\epsilon,\epsilon]\to\RR^3$, $\epsilon>0$, satisfying the conditions
from \Cref{thm:Boehm-Prautysch-Paluszny}. The idea is to construct $\bfm{\gamma}$ as a spline
of two curves $\bfm{\gamma}_1:[-\epsilon,0]\to\bfm{p}_1(\Delta)$ and 
$\bfm{\gamma}_2:[0,\epsilon]\to\bfm{p}_2(\Delta)$, such that
$\bfm{\gamma}_{|[-\epsilon,0]}=\bfm{\gamma}_1$, $\bfm{\gamma}_{|[0,\epsilon]}=\bfm{\gamma}_2$ and
$\bfm{\gamma}$ is $G^k$ continuous on $[-\epsilon,\epsilon]$. Quite clearly, 
$$
\bfm{\gamma}_i=\bfm{p}_i(\varphi_i,\psi_i),\quad i=1,2,
$$
where $(\varphi_1,\psi_1)\colon[-\epsilon,0] \to \Delta$ and $(\varphi_2,\psi_2)\colon[0,\epsilon] \to \Delta$.
Since $\bfm{\gamma}_i(0)=\bfm{c}(\tau)$, one must have $\varphi_i(0)=0$, $i=1,2$. Furthermore, 
$$
\bfm{\gamma}_i'(0)=[\varphi_i'(0),\psi_i'(0)]\nabla\bfm{p}_i(0,\psi_i(0))=
\varphi_i'(0)\frac{\partial\bfm{p}_i}{\partial u}\left(0,\psi_i(0)\right)
+\psi_i'(0)\bfm{c}'(\tau)\nparallel \bfm{c}'(\tau)
$$
implies $\varphi_i'(0)\neq 0$, $i=1,2$. Consequently, $\bfm{\gamma}_i$ can be reparameterized 
as 
$$\bfm{\gamma}_i(t)
  =\bfm{p}_i((-1)^{i}t,\tau+\sum_{j=1}^k\gamma_{ij}t^j+\widetilde{\psi}_i(t) t^{j+1}),\quad
  i=1,2, \quad t\in[0,\epsilon].
$$
Since $\bfm{\gamma}$ is a spline curve composed of $\bfm{\gamma}_1$ and $\bfm{\gamma}_2$, and the later two are
analytic, we only have to confirm the $G^k$ continuity at $t=0$. It is well know 
(see, e.g., \cite{Roman-Faa-di-Bruno-1980}) that the curve $\bfm{\gamma}$ 
is $G^k$ continuous at $t=0$ if and only if there exist
 $k$ scalar parameters $\alpha_i$, $i=1,2,\dots,k$, $\alpha_1>0$, and a lower-triangular matrix
 $M_k:=(m_{i,j}(\alpha_1,\alpha_2,\dots,\alpha_k))_{i,j=1}^k$, $m_{ii}=\alpha_1^i$, $i=1,2,\dots,k$,
 such that
 \begin{equation}\label{eq:M}
   \left(\bfm{\gamma}_2^{(j)}(0)^T\right)_{j=1}^k
   =M_k \left(\bfm{\gamma}_1^{(j)}(0)^T\right)_{j=1}^k.
 \end{equation}
 Thus we have to find a solution of the system of $3k$ nonlinear equations \eqref{eq:M}
 for $3k$ unknowns $\gamma_{ij}$, $i=1,2$, $j=1,2,\dots,k$ and $\alpha_j$, $i=1,2,\dots,k$
 with $\alpha_1>0$.\\
 In particular, matrices $M_k$, for $k=1,2,3$, are
 $$
   M_1=\left[
   \begin{array}{c}
     \alpha_1
   \end{array}
   \right], \quad
   M_2=\left[
   \begin{array}{cc}
     \alpha_1 & 0 \\
     \alpha_2 & \alpha_1^2
   \end{array}
   \right], \quad
   M_3=\left[
   \begin{array}{ccc}
     \alpha_1 & 0 & 0 \\
     \alpha_2 & \alpha_1^2 & 0\\
     \alpha_3 & 3\alpha_1\alpha_2 &\alpha_1^3
   \end{array}
   \right].
 $$

Let us now continue by considering the $G^2$ continuity of a triangular parametric cubic spline patch.
It is enough to check the $G^2$ continuity of a spline patch formed by two neighbouring
triangular parametric cubic patches $\bfm{p}_1$ and $\bfm{p}_2=R\,\bfm{p}_1$ given by the parameters \eqref{eq:best_cubic_G1}
where
\begin{align*}
R=\begin{bmatrix}
 \frac{4-5 c^2}{3 c^2-4} & 0 & \frac{4 c \sqrt{1-c^2}}{3 c^2-4} \\
 0 & 1 & 0 \\
 \frac{4 c \sqrt{1-c^2}}{3 c^2-4} & 0 & \frac{4-5 c^2}{4-3 c^2} \\
\end{bmatrix}
\end{align*}
is the reflection over the plane defined by  $\bfm{b}_{0,3,0}$, $\bfm{b}_{0,0,3}$ and the origin.
The common boundary curve $\bfm{c}$ of $\bfm{p}_1$ and $\bfm{p}_2$ is given by
$\bfm{c}:[0,1]\to\RR^3$, $\bfm{c}(v)=\bfm{p}_1(0,v)=\bfm{p}_2(0,v)$. Now we use the results from the
\Cref{subsec:Gk}.
Let us choose an arbitrary $v\in[0,1]$, small enough $\epsilon>0$ and define the curve 
$\bfm{\gamma}_v:[-\epsilon,\epsilon]\to\RR^3$ by 
\begin{align*}
  \bfm{\gamma}_v(t)=
  \begin{cases}
    \bfm{p}_1(-t,v),& t\leq 0,\\
    \bfm{p}_2\left(t+\frac{6 c^2 \left(4-5 c^2+6 c^4 (1-v) v\right) t^2}{\left(4-3 c^2\right)
   \left(4-4 c^2+3 c^4 (1-v) v\right)},v+\frac{\left(6 c^2 (1-v)-4\right) t}{4-3 c^2}+\frac{6 c^4 \left(2+9 c^4 (1-v)^2 v-3
   c^2 \left(1+v-2 v^2\right)\right) t^2}{\left(4-3 c^2\right)^2 \left(4-4 c^2+3 c^4
   (1-v) v\right)}\right),& t\ge 0.
  \end{cases}
\end{align*}
The curve $\bfm{\gamma}_v$ is regular and obviously continuous but it can be easily checked that 
it is actually $G^2$ continuous. 
By \Cref{thm:Boehm-Prautysch-Paluszny} the patches $\bfm{p}_1$ and $\bfm{p}_2$ form the $G^2$ spline patch.

\section{Quartic $G^1$ approximation}\label{section:quarticG1}  

It is clear that there are no $G^3$ continuous triangular parametric cubic spline patches approximating
(a part of) the unit sphere since they would have
to be globally polynomial. Thus we shall rise the degree of the triangular parametric patch and consider
quartic ones. In this case we have 15 control points. They are defined by \eqref{eq:vertices_of_patch}, 
by symmetries from the dihedral group $D_6$, 
and by
\begin{align*}
    \bfm{b}_{3,1,0}=\alpha\ \bfm{b}_{4,0,0}+\beta\ \bfm{b}_{0,4,0},\quad 
    \bfm{b}_{2,2,0}=\gamma\ (\bfm{b}_{4,0,0}+\bfm{b}_{0,4,0}),\quad
    \bfm{b}_{2,1,1}=\zeta\ \bfm{b}_{4,0,0}+\xi( \bfm{b}_{0,4,0}+\bfm{b}_{0,0,4}).
\end{align*}
In order to obtain at least $G^1$ continuous spline patch 
$\bfm{p}=\bfm{p}(\cdot,\cdot,\alpha,\beta,\gamma,\zeta,\xi)$, we again observe that the normal of the tangent plane of the
patch at any vertex $\bfm{b}_{4,0,0}$, $\bfm{b}_{0,4,0}$ or $\bfm{b}_{0,0,4}$ must coincide with the normal of the tangent plane 
of the sphere, i.e., it must be parallel to the radius vector of the corresponding vertex.
This implies $\alpha=\tfrac 1 2(2-2\beta+3 c^2\beta)$. Moreover, the normal of the tangent plane of the patch on the boundary curve
determined by $\bfm{b}_{4,0,0}$ and $\bfm{b}_{0,4,0}$ must lie in the plane defined by the points $\bfm{b}_{4,0,0}$, $\bfm{b}_{0,4,0}$ and 
the origin. This further implies four possibilities for parameters $(\alpha, \beta, \gamma, \zeta, \xi)$ and two of these 
possibilities again induce singular patches. The remaining two nonsingular patches are given by the two sets of parameters 
$(\alpha_i,\beta_i,\zeta_i,\xi_i)$, $i=1,2$, depending on the parameter $\gamma$ as
\begin{align}
  \alpha_1&=\frac{(6c^2-4)\gamma+c^2+2 }{4 c^2}, 
  &\alpha_2&=\frac{(9c^4-18c^2+8)\gamma-3c^4+10c^2-4}{c^2 \left(4-3 c^2\right)},\nonumber \\
  \beta_1&=\frac{2 \gamma-1}{2 c^2},
  &\beta_2&=\frac{(6c^2-8)\gamma+4}{c^2(4-3 c^2)},\nonumber \\
  \zeta_1&=\frac{-4(9c^4-18c^2+8)\gamma+9c^6-24c^4-4c^2+16}{12 c^2 \left(3 c^4-7 c^2+4\right)},
  &\zeta_2&=\frac{(9c^4-18c^2+8)\gamma-3c^4+10c^2-4}{6 c^2 \left(c^2-1\right)},\nonumber \\
  \xi_1&=\frac{-2(3c^4-3c^2-2)\gamma-c^2-2}{12 c^2 \left(1-c^2\right)},
  &\xi_2&=\frac{(9c^6-30c^4+36c^2-16)\gamma+3c^4-8c^2+8}{6 c^2 \left(3c^4-7 c^2+4\right)}.\label{eq:parameters_quartic_G1}
\end{align}
Consequently, the error functions
\begin{equation*}
  f_i(u,v,\gamma):=f(u,v,\alpha_i,\beta_i,\gamma,\zeta_i,\xi_i),\ g_i(u,v,\gamma):=g(u,v,\alpha_i,\beta_i,\gamma,\zeta_i,\xi_i),\quad i=1,2,
\end{equation*}
have to be analyzed, where
\begin{equation*}
    f(\cdot,\cdot,\alpha,\beta,\gamma,\zeta,\xi)=\|\bfm{p}(\cdot,\cdot,\alpha,\beta,\gamma,\zeta,\xi)\|_2^2-1, \quad
    g(\cdot,\cdot,\alpha,\beta,\gamma,\zeta,\xi)=\sqrt{f(\cdot,\cdot,\alpha,\beta,\gamma,\zeta,\xi)+1}-1.
\end{equation*}
Let us consider three particular cases arising from polyhedra with parameters $c$ from \eqref{eq:particular_c}.

\subsection{Tetrahedron} 
In this case $c=\tfrac{2\sqrt{2}}{3}$. Consider $f_1$ first. Take an arbitrary point $(u,v)\in\Delta$.
The function $f_1(u,v,\cdot)$ is a quadratic function with a positive leading coefficient and its
minimum located on $(-\infty,\tfrac{1}{2}]$. 
Therefore $f_1(u,v,\cdot)$ is an increasing function on $[\tfrac 1 2,\infty)$. Furthermore, the value
$\gamma_{f_1}=\tfrac{2189+108\sqrt{2291}}{7602}\approx 0.967950$ is the only solution of the equation 
$f_1\left(\tfrac 1 3,\tfrac 1 3,\gamma\right)=-f_1\left(\tfrac 1 2,\tfrac 1 2,\gamma\right)$ on $[0,\infty)$. 
Using \Cref{lem:ftog} one can confirm that
$f_1\left(\tfrac 1 2,\tfrac 1 2,\gamma_{f_1}\right)$ and $f_1\left(\tfrac 1 3,\tfrac 1 3,\gamma_{f_1}\right)$ 
are global minimum and maximum
of $f_1(\cdot,\cdot,\gamma_{f_1})$ on $\Delta$.\\
If $\gamma>\gamma_{f_1}$ then $f_1\left(\tfrac 1 3,\tfrac 1 3,\gamma\right)>f_1\left(\tfrac 1 3,\tfrac 1 3,\gamma_{f_1}\right)$ 
since $f_1\left(\tfrac 1 3,\tfrac 1 3,\cdot\right)$ is an increasing function. On the other hand, if $\gamma<\gamma_{f_1}$,
then again by the monotonicity we have $f_1\left(\tfrac 1 2,\tfrac 1 2,\gamma\right)<f_1\left(\tfrac 1 2,\tfrac 1 2,\gamma_{f_1}\right)$.
This confirms that $\bfm{p}(\cdot,\cdot,\gamma_{f_1})$ is an optimal approximant.\\
%An easy evaluation reveals that
%$f_1\left(\tfrac 1 3 ,\tfrac 1 3,\gamma^*\right)=\tfrac{64(12496-243\sqrt{2291})}{1605289}\approx 0.034841$.
Let us now consider the function $f_2$. We shall see that the approximant arising from the
second set of parameters (indexed by $2$ in \eqref{eq:parameters_quartic_G1})
is inferior to the approximant arising from the first set of parameters 
(indexed by $1$ in  \eqref{eq:parameters_quartic_G1})
for any parameter $\gamma$. 
If $\gamma>\widetilde{\gamma}$, where $\widetilde{\gamma}=\tfrac{22229-216\sqrt{2291}}{7602}\approx 1.56410$, 
then $f_2(\tfrac 1 2,\tfrac 1 2,\gamma)<f_1(\tfrac 1 2,\tfrac 1 2,\gamma_{f_1})$.
If $\gamma<\widetilde{\gamma}$ then $f_2(\tfrac 1 3,\tfrac 1 3,\gamma)>f_1(\tfrac 1 3,\tfrac 1 3,\gamma_{f_1})$.
If $\gamma=\widetilde{\gamma}$ then 
$f_2(\tfrac 1 2-\tfrac 1 2\sqrt{\tfrac{27\sqrt{2291}-979}{2649}},0,\gamma)\approx
f_2(0.328037,0,\gamma)
<f_1(\tfrac 1 2,\tfrac 1 2,\gamma_{f_1})$ and the desired conclusion follows.\\
An almost identical analysis can be done for $g_1$ and the resulting
optimal parameter is $\gamma_{g_1}=\frac{9587-2916 \sqrt{3}}{4686}\approx 0.968062$.

\subsection{Octahedron} 
If the underlying polyhedron is an octahedron, we have $c=\tfrac{\sqrt{6}}{3}$. Let us again first analyze the 
error function $f_1$.\\
Observe that functions $f_1(\tfrac 1 2, \tfrac 1 2,\cdot)$ and $f_1(\tfrac 1 3, \tfrac 1 3,\cdot)$ are both quadratic again. 
The first one is increasing on $[-\tfrac 1 6,\infty)$ and the second one increases on $[\tfrac{1}{8},\infty)$. Moreover, it is easy to check that
$\gamma_{f_1}=\tfrac{47+36\sqrt{974}}{1510}\approx 0.775181$ is the only solution of the equation 
$f_1(\tfrac 1 3,\tfrac 1 3,\gamma)=-f_1(\tfrac 1 2,\tfrac 1 2,\gamma)$ on $[0,\infty)$. 
The extrema of $f_1$ can be quite easily obtained by using \Cref{lem:ftog} again.
We get
$f_1(\tfrac 1 2,\tfrac 1 2,\gamma_{f_1})\le f_1(u,v,\gamma_{f_1})\le f_1(\tfrac 1 3 ,\tfrac 1 3,\gamma_{f_1})$
for any $(u,v)\in \Delta$.\\
If $\gamma>\gamma_{f_1}$ then $f_1(\tfrac 1 3,\tfrac 1 3,\gamma)>f_1(\tfrac 1 3,\tfrac 1 3,\gamma_{f_1})$ and
if $\gamma<\gamma_{f_1}$ then $f_1(\tfrac 1 2,\tfrac 1 2,\gamma)<f_1(\tfrac 1 3,\tfrac 1 3,\gamma_{f_1})$
since $f_1(\tfrac 1 3,\tfrac 1 3,\cdot)$ and $f_1(\tfrac 1 2,\tfrac 1 2,\cdot)$ are both
increasing functions, so $\bfm{p}(\cdot,\cdot,\gamma_{f_1})$ is the best
approximant.\\
Similarly as in the case of tetrahedron, we can now show that the second set of parameters in
\eqref{eq:parameters_quartic_G1} induces worse approximant for any $\gamma$. 
%Namely,  
%if $\gamma>\widetilde{\gamma}$, where $\widetilde{\gamma}=\tfrac{9(409-8\sqrt{974})}{1510}$ 
%then $f_2(\tfrac 1 2,\tfrac 1 2,\gamma)<f_1(\tfrac 1 2,\tfrac 1 2,\gamma_{f_1})$.
%If $\gamma<\widetilde{\gamma}$ then $f_2(\tfrac 1 3,\tfrac 1 3,\gamma)>f_1(\tfrac 1 3,\tfrac 1 3,\gamma_{f_1})$ 
%and if $\gamma=\widetilde{\gamma}$ then 
%$f_2(\tfrac 1 2-\sqrt{\tfrac{27\sqrt{974}-619}{3666}},0,\gamma)<f_1(\tfrac 1 2,\tfrac 1 2,\gamma_{f_1})$.
If we perform a similar analysis for $g_1$, the optimal parameter turns out to be\\
$\gamma_{g_1}=\frac{1}{538} \left(209+768 \sqrt{3}-12 \sqrt{6126+1512 \sqrt{3}}\right)
\approx 0.775181$.

\begin{figure}[!htb]
  \minipage{0.33\textwidth}
    \centering
    \includegraphics[width=1.0 \linewidth]{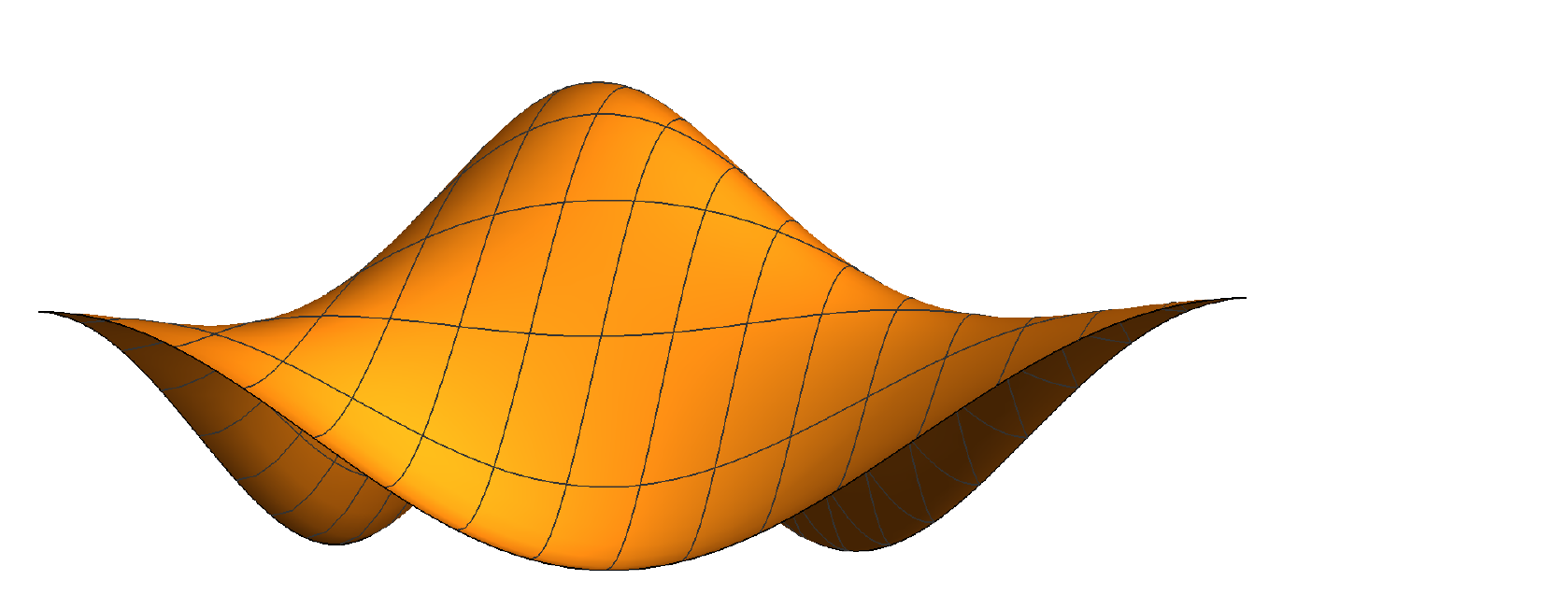}
  \endminipage\hfill
  \minipage{0.33\textwidth}
    \centering
    \includegraphics[width=1.0 \linewidth]{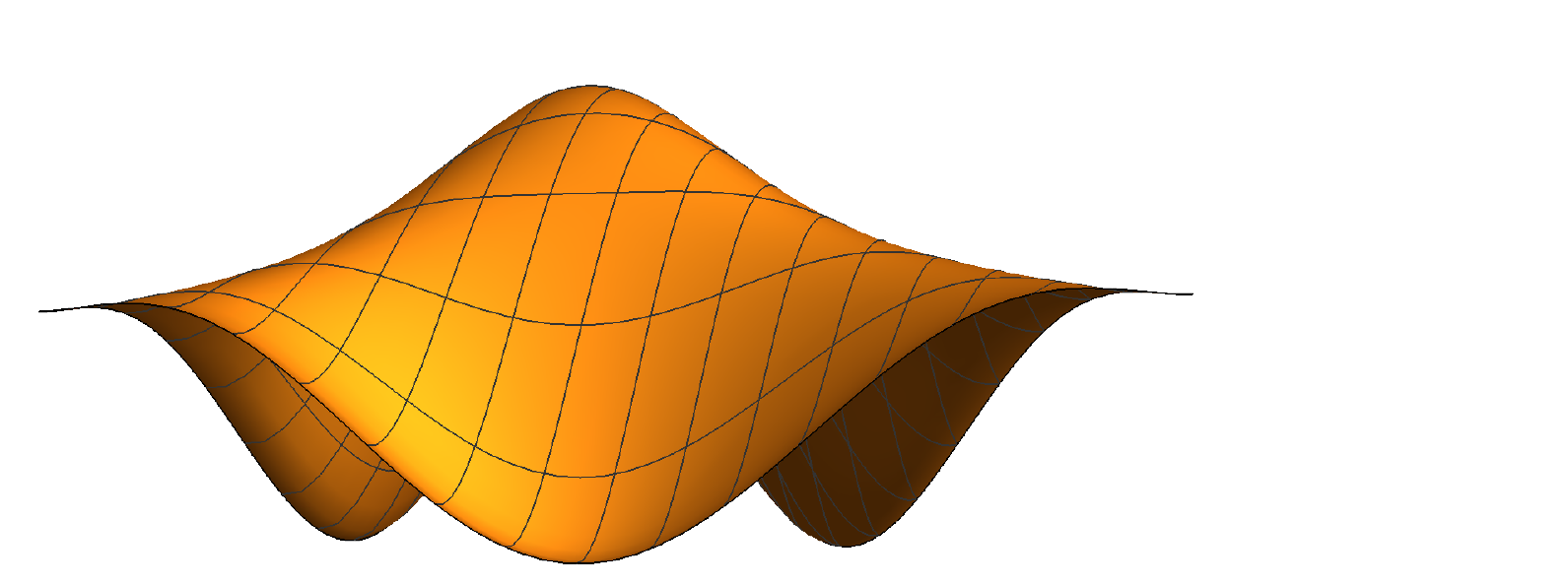}
  \endminipage\hfill
  \minipage{0.33\textwidth}
    \centering
    \includegraphics[width=1.0 \linewidth]{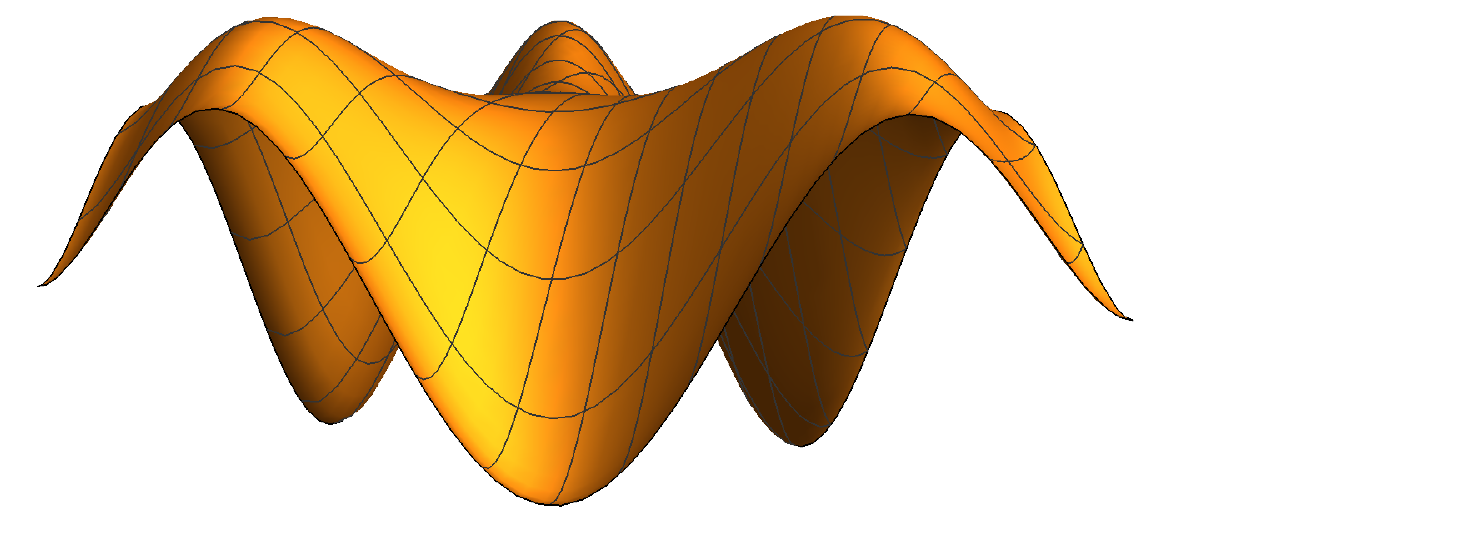}
  \endminipage\hfill
  \caption{The graphs of simplified error functions $f$ of the best quartic interpolants over one triangle of the tetrahedron, the octahedron 
  and the icosahedron triangulation, respectively.} 
  \label{fig:G1_quartic_errors}
  %with maximal errors $3.4\times 10^{-2}$, $2.0\times 10^{-3}$ and $6.3\times10^{-5}$, respectively.}
\end{figure}

\subsection{Icosahedron}

If the underlying polyhedron is an icosahedron, then by \eqref{eq:particular_c} we have $c=\sqrt{\tfrac{2(5-\sqrt{5})}{15}}$.
The analysis of the error function is in this case more complicated then in the previous ones. 
The reason is that the optimal interpolant induces the error function which does not have a maximum 
at the barycentre of the triangle (see \Cref{fig:G1_quartic_errors}). More precisely, the positive solution of the quadratic equation 
$f_1(\tfrac 1 3,\tfrac 1 3,\gamma)=-f_1(\tfrac 1 2,\tfrac 1 2,\gamma)$ is 
$$
\gamma=\frac{-1931-1100 \sqrt{5}+18 \sqrt{6 \left(556615+248877 \sqrt{5}\right)}}{6 \left(5771+2508 \sqrt{5}\right)}:=\gamma_0\approx 0.617027,
$$
and it induces the interpolant for which $(\tfrac 1 3,\tfrac 1 3)\in\Delta$ does not provide a maximum 
but $(\tfrac{1}{2},\tfrac{1}{2})\in\Delta$ implies a minimum.
Furthermore, the point where $f_1(\cdot,\cdot,\gamma)$ attains the maximal value depends on the parameter $\gamma$.

It is easy to see that $\gamma_{f_1}\ge \tfrac 1 2$, since $|f_1(\tfrac 1 2,\tfrac 1 2,\gamma)|$ is too big for 
$\gamma<\tfrac 1 2$ and consequently the error function can be easily improved.
For every $(u,v)\in\Delta$ the function $f_1(u,v,\cdot)$ is a quadratic polynomial 
and it can be shown that it is increasing for $\gamma\ge \tfrac 1 2$. Consequently, 
$\gamma_{f_1}\in [\tfrac 1 2, \gamma_0]$ and the parameter $\gamma_{f_1}$ can be found using the bisection method as follows: 
For $\gamma$ being the midpoint of the interval $[\tfrac 1 2, \gamma_0]$ we easily calculate the maximum $M_\gamma$ 
and the minimum $m_\gamma$ of the function $f_1(\cdot,\cdot,\gamma)$ on $\Delta$. 
If $M_\gamma>|m_\gamma|$ then $\gamma_{f_1}\in [\tfrac 1 2,\gamma]$ else $\gamma_{f_1}\in [\gamma,\gamma_0]$ 
and we can continue with the procedure.\\
The value $\gamma_{f_1}$ can also be computed directly. By using \Cref{lem:ftog} again,
we can prove that for every $\gamma$ the extrema of $f_1(\cdot,\cdot,\gamma)$ on $\Delta$ appear on the boundary 
of $\Delta$ or on its medians. 
%($u=0$, $v=0$ or $v+v=1$) or on the medial of the triangle $\Delta$ ($u=v$, $v=1-2u$ or $u=1-2v$).
It is straightforward to observe that the solutions of the equations
$f_1(\tfrac 1 3,\tfrac 1 3,\gamma_1)=0$ and $f_1(\tfrac 1 2,\tfrac 1 2,\gamma_2)=0$ are 
$$
\gamma_1=\frac{1}{96} \left(29-13 \sqrt{5}+9 \sqrt{6 \left(5+\sqrt{5}\right)}\right)\approx 0.617012,
\gamma_2=\frac{1}{6} \left(-1+2 \sqrt{10-2 \sqrt{5}}\right)
\approx 0.617047.
$$
We shall see now that $\gamma_{f_1}\in [\gamma_1,\gamma_2]$. For every $\gamma\in [\gamma_1,\gamma_2]$ the function 
$f_1(\cdot,\cdot,\gamma)$ has the minimum on $\Delta$ at $(\tfrac 1 2,0)$. To show this we define the function $k(u,v,t)=f_1(u,v,(\gamma_2-\gamma_1)t+\gamma_1)$ and the inequality $f_1(u,v,\gamma)-f_1(\tfrac 1 2,0,\gamma)\ge 0$ for all $(u,v)\in\Delta$ and all $\gamma\in[\gamma_1,\gamma_2]$ is equivalent to the inequality $k(u,v,t)-k(\tfrac 1 2,0,t)\ge 0$ for all $(u,v)\in\Delta$ and all $t\in[0,1]$. The later is true since
$k(u,v,t)-k(\tfrac 1 2,0,t)=a_0(u,v)(1-t)+a_1(u,v)t(1-t)+a_2(u,v)t^2$
and polynomials $a_0$, $a_1$ and $a_2$ are non-negative on $\Delta$.\\
We now know that the minimum of $f_1(\cdot,\cdot,\gamma)$ is always at $(\tfrac 1 2,0)$.
On the other hand, it is by \Cref{lem:ftog} enough to consider one of the subtriangles
on \Cref{fig:delta_omega} (on the left) as a subdomain. Again by \Cref{lem:ftog} there is no 
extrema in the interior of this domain.
Consequently, the maximum on the subdomain must appear for $u=v$, $v=0$ or $v=1-2u$. 
Some analysis leads to the system of polynomial (in)equations 
\begin{equation*}
f_1(u,u,\gamma)=-f_1(\tfrac 1 2,0,\gamma),\quad 
%& f_1(u,0,\gamma)&=-f_1(\tfrac 1 2,0,\gamma),\\
\frac{\partial}{\partial u}f_1(u,u,\gamma)=0,\quad
%&\frac{\partial}{\partial u}f_1(u,0,\gamma)&=0,\\
\frac{\partial^2}{\partial u^2}f_1(u,u,\gamma)<0. 
%&\frac{\partial^2}{\partial u^2}f_1(u,0,\gamma)&<0.
\end{equation*}
The only admissible solution is $\gamma_{f_1}\approx 0.617022$ and $u^*=v^*\approx 0.139979$.\\
As in the cases of thetrahedron and octahedron, the optimal parameter $\gamma_{g_1}$ can be found similarly as $\gamma_{f_1}$. Unfortunately, the analytical approach explained above is not successful for $g_1$,
but one can still use the method of bisection explained previously and the numerical approximation $\gamma_{g_1}\approx 0.617022$ follows.\\
Quite straightforward but technically challenging computations show that the functions $f_2$ and $g_2$
provide parameters $\gamma_{f_2}$ and $\gamma_{g_2}$ which imply inferior approximants but the details will be
omitted here.

\begin{figure}[!htb]
  \minipage{0.33\textwidth}
    \centering
    \includegraphics[width=0.75 \linewidth]{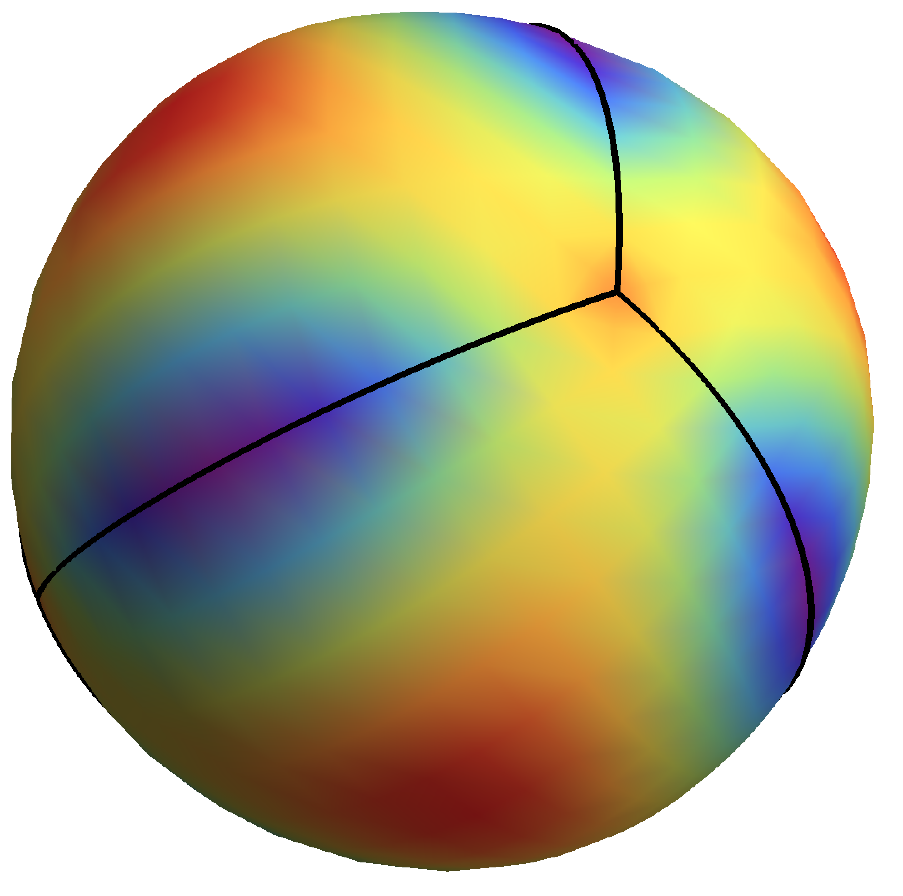}
  \endminipage\hfill
  \minipage{0.33\textwidth}
    \centering
    \includegraphics[width=0.75\linewidth]{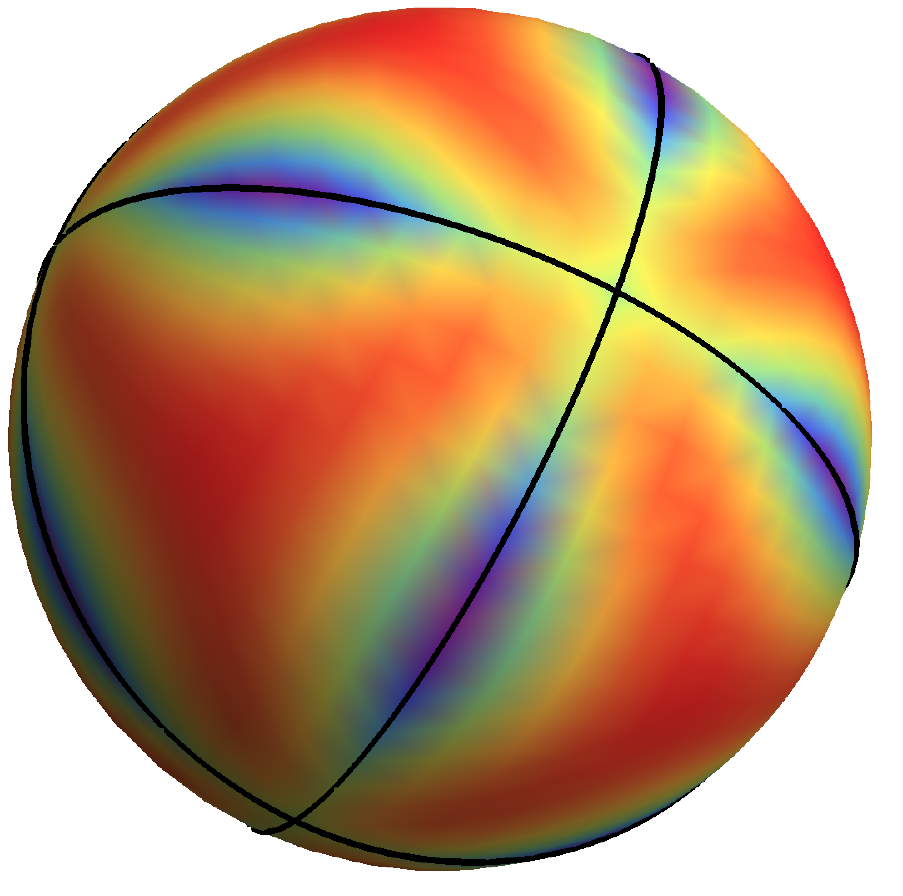}
  \endminipage\hfill
  \minipage{0.33\textwidth}
    \centering
    \includegraphics[width=0.75\linewidth]{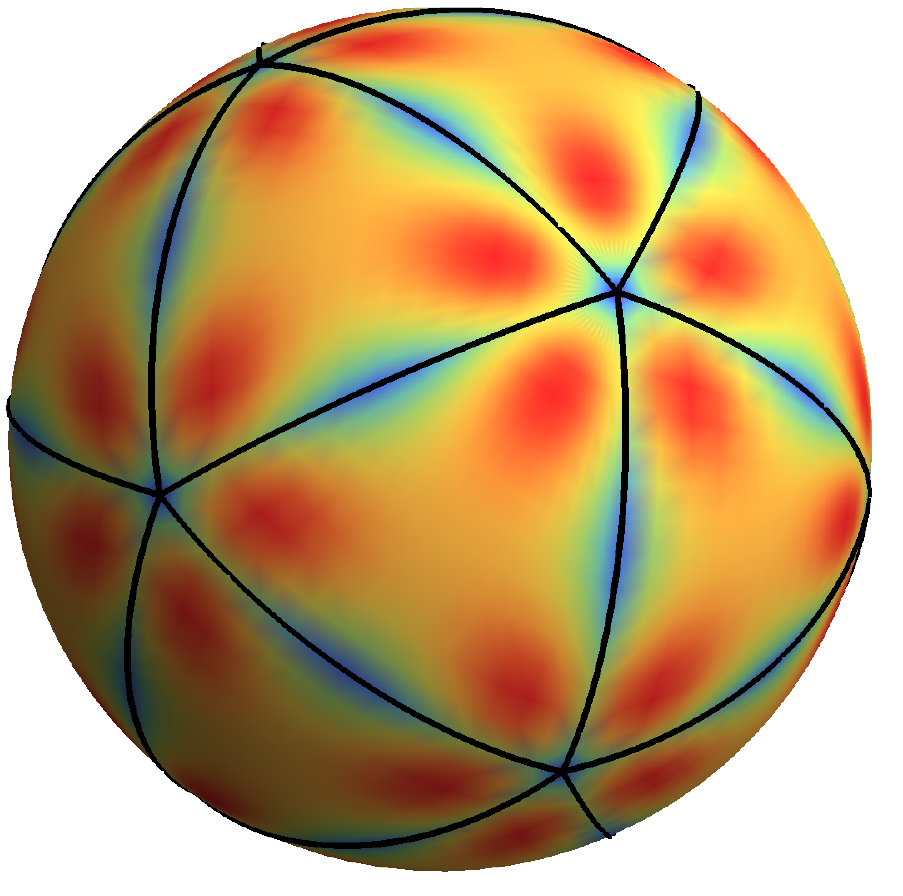}
  \endminipage\hfill
  
  \caption{Figures of the best quartic $G^1$ approximants of the whole sphere 
  based on the underlying tetrahedron, octahedron and icosahedron (from left to right)
  together with the Gaussian curvatures (red regions indicate higer curvature). 
  \label{fig:best_quart_G1}}
\end{figure}

The results for all three cases are collected in \Cref{tab:G14}.
The optimal $G^1$ quartic spline patches together with their Gaussian curvatures are shown 
on \Cref{fig:best_quart_G1}.

\renewcommand{\arraystretch}{1.5}
\begin{table}[h]
  \begin{equation*}
    \begin{array}{|l|r|r|r|r|r|r|r|r|}\hline
        & \multicolumn{1}{|c|}{\alpha_{1}} & \multicolumn{1}{|c|}{\beta_{1}} & \multicolumn{1}{|c|}{\gamma_{1}} 
        & \multicolumn{1}{|c|}{\zeta_{1}}  & \multicolumn{1}{|c|}{\xi_{1}} 
        & \multicolumn{1}{|c|}{d_r} &\multicolumn{1}{|c|}{K_{min}} & \multicolumn{1}{|c|}{K_{max}}\\ \hline
        \text{tetrahedron} 
        & 1.175523 & 0.526570 & 0.968062 
        & 2.053140 & 1.313741
        & 0.017296 &0.68 &1.24\\ \hline
        \text{octahedron}
        & 1.000000 & 0.412772 &0.775181
        & 1.000000 & 0.550362
        & 0.001019 & 0.93 & 1.03\\ \hline
        \text{icosahedron} 
        & 0.857991 & 0.317543 &0.617022
        & 0.659094 & 0.344164
        & 0.000017 & 0.99 & 1.00\\ \hline
    \end{array}
  \end{equation*}
  \caption{Optimal parameters $\alpha_1$, $\beta_1$, $\gamma_1$, $\zeta_1$ and $\xi_1$ 
  according to the radial error, radial distances $d_r$ and the corresponding
  minimal and maximal Gaussian curvatures $K_{min}$, $K_{max}$ for the optimal triangular parametric
  $G^1$ quartic patches. The parameters according to the simplified radial error are almost the same
  and they are omitted in this table.} \label{tab:G14}
\end{table}

\section{Quartic $G^2$ approximation}\label{section:quarticG2}

Similarly as in the cubic case, we will prove that each $G^1$ quartic spline approximant derived in the previous
section is actually a $G^2$ spline. A similar approach will be used as in the case of $G^1$ cubic spline patches
and it relies on the results of \Cref{thm:Boehm-Prautysch-Paluszny} again.\\
Let $\bfm{p}_1(\cdot,\cdot,\gamma):=\bfm{p}(\cdot,\cdot,\alpha_1,\beta_1,\gamma,\zeta_1,\xi_1)$ 
be the quartic $G^1$ interpolant determined by the parameters in \eqref{eq:parameters_quartic_G1}. 
Furthermore, let 
\begin{align*}
  R=\begin{bmatrix}
        \frac{4-5 c^2}{3 c^2-4} & 0 & \frac{4 c \sqrt{1-c^2}}{3 c^2-4} \\
        0 & 1 & 0 \\
        \frac{4 c \sqrt{1-c^2}}{3 c^2-4} & 0 & \frac{4-5 c^2}{4-3 c^2} \\
    \end{bmatrix}
\end{align*}
be the matrix representing the reflection over the plane defined by the coordinate origin, $\bfm{b}_{0,4,0}$ and $\bfm{b}_{0,0,4}$ and let $\bfm{p}_2(\cdot,\cdot,\gamma)=R\,\bfm{p}_1(\cdot,\cdot,\gamma)$.
Since $\bfm{p}_1(0,v,\gamma)=R\,\bfm{p}_1(0,v,\gamma)$, the spline patch composed by $\bfm{p}_1$ and $\bfm{p}_2$
is at least continuous. Let the curve $\bfm{\gamma}_v$, $v\in [0,1]$, be defined as
\begin{align*}
\bfm{\gamma}_v(t)=
\begin{cases}
\bfm{p}_1(-t,v,\gamma),& t\le 0,\\
\bfm{p}_2(t,v+\frac{2 t \left(2-3 c^2+3 c^2 v\right)}{3 c^2-4}+\frac{2 t^2 \left(27 c^6 (1-v)^2 v-16-6 c^2 (4 v-7)+9 c^4 \left(2 v^2+v-3\right)\right)
   \left(2+c^2-4 \gamma \right)}{\left(4-3 c^2\right)^2 \left(6 c^4 (v-1) v \gamma +2 \left(2 v^2-2 v+1\right) (2 \gamma -1)+c^2 \left(2+v^2
   (7-18 \gamma )-4 \gamma +v (18 \gamma -7)\right)\right)},\gamma),& t\ge 0.
\end{cases}
\end{align*}
For $\gamma>\tfrac{6}{10}$ the curve $\bfm{\gamma}_v$ is well defined for all $v\in[0,1]$. Note that
this includes all three optimal parameters $\gamma$ for the tetrahedron, octahedron and icosahedron case
derived in the previous section.
Clearly $\bfm{\gamma}_v(0)=\bfm{p}_1(0,v,\gamma)$ 
and it can be verified that $\bfm{\gamma}_v$ is a $G^2$ curve. Therefore, 
by \Cref{thm:Boehm-Prautysch-Paluszny} the patches $\bfm{p}_1$ and $\bfm{p}_2$ form the $G^2$ spline patch.
This implies that the best $G^1$ quartic spline approximant from the previous section is actually 
the best $G^2$ quartic spline approximant.

\section{Closure}\label{sec:closure}
In surface design it is important to develop efficient and accurate algorithms for parametric polynomial 
approximation of (a part of) a sphere. The quality relies on measuring the distance between parametric surfaces 
which is not a trivial task. In this study a similar approach as in the case of parametric curves
was used and the distance between parametric surfaces was measured as a (simplified) radial error. 
The optimal approximation of a canonical surface, i.e., an equilateral spherical triangle was studied 
and the results were obtained for low degree parametric polynomial patches with a particular order of 
geometric continuity. In case when equilateral spherical triangles can be put
together to form the unit sphere, optimal parametric polynomial spline patches were studied. 
For almost all of them closed form solutions were obtained which makes them useful in practical applications.\\
For the future work the study of optimal approximation of isosceles spherical triangles or even general 
ones is planned. This would lead to the optimal approximation of some other canonical surfaces, 
such as ellipsoide or hyperboloid. One could also study an optimal approximation by tensor product parametric
polynomial patches. The first step in this direction was done in \cite{Eisele-1994-best-biquadratic}.\\

{\noindent \sl Acknowledgments.}
  The first author was supported by the Slovenian Research Agency
  program P1-0292 and the grants J1-8131, N1-0064, and N1-0083. The second author was supported in part by the program
  P1-0288 and the grant J1-9104 by the same agency.

\bibliography{bib_cagd}
\end{document}